\documentclass[11pt,fleqn,reqno]{article}
\usepackage{epsfig,amssymb,latexsym,amsmath,pifont,multicol,mathtools,amsmath,verbatim,amsthm,float,lipsum}
\usepackage{caption} 

\usepackage[utf8]{inputenc}
\usepackage[T1]{fontenc}

\usepackage[varg]{txfonts}
\let\mathbb=\varmathbb

\usepackage[sans]{dsfont}

\usepackage[left=2.5cm, right=2.5cm, top=2.5cm, bottom=2.5cm]{geometry}

\usepackage[%
cal=euler,
bb=fourier,
scr=zapfc,
]
{mathalfa}

\usepackage[font=small,labelfont=bf]{caption}
\usepackage{subcaption}
\usepackage{graphicx}
\graphicspath{{Figures/}} 
\usepackage{acronym}
\usepackage{latexsym}
\usepackage{paralist}
\usepackage{wasysym}
\usepackage{xspace}
\usepackage{framed}
\usepackage{palatino,pxfonts}
\usepackage{authblk}
\usepackage{enumitem}
\usepackage[numbers,sort&compress]{natbib}

\usepackage[dvipsnames,svgnames]{xcolor}
\colorlet{MyBlue}{DodgerBlue!75!Black}
\colorlet{MyGreen}{DarkGreen!95!Black}

\newcommand{\MS}[1]{{\color{red}[\textbf{MS}:\;#1]}}
\newcommand{\DL}[1]{{\color{blue}[\textbf{DL}:\; #1]}}

\usepackage{hyperref}
\hypersetup{
colorlinks=true,
linktocpage=true,
pdfstartview=FitH,
breaklinks=true,
pdfpagemode=UseNone,
pageanchor=true,
pdfpagemode=UseOutlines,
plainpages=false,
bookmarksnumbered,
bookmarksopen=false,
bookmarksopenlevel=1,
hypertexnames=true,
pdfhighlight=/O,
urlcolor=MyBlue!60!black,linkcolor=MyBlue!70!black,citecolor=DarkGreen!70!black, 
pdftitle={},
pdfauthor={},
pdfsubject={},
pdfkeywords={},
pdfcreator={pdfLaTeX},
pdfproducer={LaTeX with hyperref}
}

\numberwithin{equation}{section}  
\usepackage[sort&compress,capitalize,nameinlink]{cleveref}
\crefname{example}{Ex.}{Exs.}

\crefrangeformat{equation}{\upshape(#3#1#4)\textendash(#5#2#6)}

\usepackage{mathtools}
\usepackage[para,online,flushleft]{threeparttable}
\usepackage{multirow}

\usepackage{physics}
\usepackage{float}
\newcommand{\eps}{\varepsilon}
\DeclareMathOperator*{\argmin}{argmin}

\DeclareMathOperator{\diag}{diag}

\DeclareMathOperator{\dist}{dist}
\DeclareMathOperator{\sign}{sign}
\DeclareMathOperator{\dif}{d\!}

\DeclareMathOperator{\range}{ran}

\DeclareMathOperator{\cond}{cond}

\DeclareMathOperator{\Id}{Id}

\DeclareMathOperator{\prox}{\textsf{prox}}
\newcommand{\proj}{\operatorname{Proj}}
\newcommand{\Sym}{\operatorname{Sym}}

\newcommand{\Retract}{\operatorname{Retr}}
\newcommand{\bA}{{\mathbf A}}
\newcommand{\bB}{{\mathbf B}}
\newcommand{\bQ}{{\mathbf Q}}

\newcommand{\bX}{\mathbf{X}}

\newcommand{\bV}{\mathbf{V}}

\newcommand{\UU}{\mathbf{U}}
\newcommand{\bx}{\mathbf{x}}

\newcommand{\bu}{\mathbf{u}}

\newcommand{\bS}{\mathbf{S}}
\newcommand{\bZ}{\mathbf{Z}}
\newcommand{\bT}{\mathbf{T}}

\newcommand{\bW}{\mathbf{W}}
\newcommand{\bv}{\mathbf{v}}
\newcommand{\bK}{\mathbf{K}}

\renewcommand{\emptyset}{\varnothing}

\newcommand{\scrC}{\mathcal{C}}

\newcommand{\scrM}{\mathcal{M}}

\newcommand{\scrO}{\mathcal{O}}

\newcommand{\scrS}{\mathcal{S}}

\newcommand{\scrU}{\mathcal{U}}

\newcommand{\setX}{\mathsf{X}}
\newcommand{\setY}{\mathsf{Y}}

\newcommand{\R}{\mathbb{R}}

\newcommand{\N}{\mathbb{N}}

\newcommand{\bC}{{\mathbf{C}}}
\newcommand{\bD}{{\mathbf{D}}}

\newcommand{\St}{\operatorname{St}}

\newcommand{\la}{\lambda}
\newcommand{\La}{\Lambda}

\newcommand{\ga}{\gamma}

\newcommand{\al}{\alpha}

\newcommand\snorm[1]{\Vert#1\Vert}

\newcommand\sabs[1]{\left\lvert#1\right\rvert}
\newcommand\inner[2]{\langle#1,#2\rangle}

\usepackage{algorithm}
\usepackage{algpseudocode}								

\theoremstyle{plain}
\newtheorem{theorem}{Theorem}
\newtheorem{corollary}[theorem]{Corollary}
\newtheorem*{corollary*}{Corollary}
\newtheorem{lemma}[theorem]{Lemma}
\newtheorem{proposition}[theorem]{Proposition}

\newtheorem{problem}{Problem}
\theoremstyle{definition}
\newtheorem{definition}[theorem]{Definition}
\newtheorem*{definition*}{Definition}
\newtheorem*{problem*}{Problem}
\newtheorem{assumption}{Assumption}

\theoremstyle{remark}
\newtheorem{remark}{Remark}
\newtheorem*{remark*}{Remark}
\newtheorem*{notation*}{Notational remark}
\newtheorem{example}{Example}

\numberwithin{theorem}{section}
\numberwithin{remark}{section}
\numberwithin{example}{section}

\usepackage[section]{placeins}

\makeatletter
\def\fps@figure{!htbp}
\makeatother

\title{Computing Data-Aware Diagonal Frame Decompositions}
\author[1]{Andrea Ebner}
\author[2]{Dirk Lorenz}
\author[3]{Mathias Staudigl}

\affil[1]{Interdisciplinary Center for Scientific Computing, Heidelberg University, Heidelberg, Germany\\
\href{mailto:andrea.ebner@uibk.ac.at}{andrea.ebner@uibk.ac.at}}

\affil[2]{Center for Industrial Mathematics,
University of Bremen, Bremen, Germany \\
\href{d.lorenz@uni-bremen.de}{d.lorenz@uni-bremen.de}}

\affil[3]{Department of Business Analytics and Decision Sciences,
WU Vienna, Vienna, Austria \\
\href{Mathias.Staudigl@wu.ac.at}{Mathias.Staudigl@wu.ac.at}}
\date{\today}

\begin{document}

\maketitle

\begin{abstract}%
We present a method for computing a specific type of matrix decomposition that combines a structured framework with data-driven adaptation. This framework allows for an efficient and stable representation of signals or images using frames by identifying sets of vectors that naturally capture the essential features of the given data. The resulting decomposition enables the development of regularization techniques that enhance the quality of reconstruction by preserving important signal characteristics while suppressing noise. We develop an optimization-based algorithm to compute the decomposition and establish convergence properties of the resulting iterates. Numerical experiments on synthetic data and tomographic image reconstruction illustrate the behavior of the proposed method and its performance compared to conventional spectral regularization methods.
\end{abstract}

\section{Introduction}
\label{sec:Intro}
%

In this paper, we consider a matrix decomposition problem that we call \emph{diagonal tight frame decomposition}~\cite{Fr19} with data-aware adaptation for sparse representation. The decomposition is related to the singular value decomposition, with the difference that the corresponding singular vectors form a tight frame (for a precise definition, see Section~\ref{sec:frames}) instead of an orthonormal basis. The problem is formulated as follows:
\begin{problem} \label{decom_problem}
Given a matrix $\bA \in \R^{N \times M}$ and a data set $\bX = [\bx_1, \bx_2, \ldots, \bx_S] \in \R^{M \times S}$, find matrices $\UU, \bD$, and $\bV$ satisfying the following properties:
\begin{enumerate}
\item \textbf{Tightness:} $\UU$ and $\bV$ are tight frames,
\item \textbf{Scalability:} $\bD$ is diagonal,
\item \textbf{Decomposition:} $\bA^{\top} \bV = \UU \bD$,
\item \textbf{Sparsity:} the representations $\UU^{\top} \bX$ are sparse.
\end{enumerate}
\end{problem}

The motivation for this problem comes from the use of frame decompositions for the representation and inversion of linear operators in Hilbert spaces. This diagonalized structure allows the evaluation and inversion of the operator to be analysed and implemented in terms of frame operators and scalar multiplication. The sparsity property incorporates prior information about the underlying signal class and can help stabilize the representation. Such decompositions have been used, for example, in the regularization of inverse problems \cite{ebner2023convergence, ebner2023regularization, donoho1995nonlinear, Fr19,frikel2018efficient, candes2002recovering}.

If $\UU$ and $\bV$ are orthonormal bases for $\R^M$ and $\R^N$, respectively, then $\bA = \bV \bD \UU^{\top}$ 
is a singular value decomposition. However, the singular vectors generally do not provide sparse representations for a given signal class. They are determined solely by the matrix $\bA$ and therefore cannot be adapted to different signal distributions.
In contrast, tight frames provide stable and redundant representations that are well suited for sparse signal representations. Moreover, their dual frames are particularly easy to compute, since the canonical dual of a tight frame differs only by a multiplicative constant. The diagonal matrix $\bD$ plays a role analogous to the singular values and is therefore often referred to as the matrix of quasi-singular values. It represents the scaling action of the operator in the corresponding frame coordinates.
Replacing orthonormal bases by tight frames therefore provides greater flexibility in the construction of such decompositions, allowing additional properties to be imposed. In particular, our sparsity condition introduces a data-dependent criterion: rather than constructing a decomposition solely from the algebraic properties of $\bA$, we seek a frame $\UU$ whose coefficients provide sparse representations given the data $\bX$. This is motivated by the observation that suitable representations can capture the structure of signals and images more efficiently and incorporate prior information that is relevant for inverse problems \cite{Selesnick05,Chaux2006}.
We note in passing that the diagonal matrix decomposition problems has applications far beyond imaging. In fact, frames are widely used for this purpose in applied mathematics, computer science, and engineering, with applications including sampling theory \cite{eldar2003sampling}, data quantization \cite{blum2010sobolev,bodmann2007frame}, image processing \cite{candes2004new,kutyniok2011compactly}, and coding theory \cite{strohmer2003grassmannian}. 

\subsection{Applications}

\paragraph{Application 1: Inverse problems.}
In inverse problems, one seeks to recover an unknown parameter $x\in\setX$ from noisy measurements
\begin{equation}\label{eq}
\bK x + n = y^\delta.
\end{equation}
$\bK:\setX\to\setY$ is a linear operator between Hilbert spaces $\setX$ and $\setY$, $n$ denotes a bounded noise term, and $y^\delta\in\setY$ is the observed data. The goal is to recover the unknown $x$ from the noisy measurements $y^\delta$. 

The operator $\bK$ is typically not continuously invertible, making the reconstruction problem ill-posed and necessitating regularization. Given a discretization $\bA \in \mathbb{R}^{N\times M}$ of $\bK$, we follow the \emph{discretize-then-optimize} methodology from numerical analysis. Since $\bA$ is typically ill-conditioned, direct inversion is highly sensitive to measurement noise and therefore does not provide a stable reconstruction.

Given the decomposition $(\UU,\bV,\bD)$ searched for in Problem~\ref{decom_problem}, we have $\bA^{\top} \bV = \UU \bD$ for tight frames $\UU$ and $\bV$ and diagonal matrix $\bD$. Let $x^+ \in \ker(\bA)^\perp$ be the least-squares solution of $\bA x = y$ for noiseless observations $y \in \range(\bA)$. By the sparsity constraint, the frame coefficients of $x^+$ are assumed to be sparse with respect to $\UU$. In particular, given a dual frame $\bar{\UU}$ of $\UU$, the duality relation for frames gives
\[
x^+ = \bar{\UU}\UU^{\top} x^+
= \sum_{i=1}^K \inner{x^+}{u_i}\bar{u}_i,
\]
where $u_i$ and $\bar{u}_i$ denote the columns of $\UU$ and $\bar{\UU}$, respectively. Thus, if the frame coefficients $\inner{x^+}{u_i}$ are sparse, the signal can be represented accurately using only a small number of frame elements.
The decomposition further yields the relation $\bD^\dagger \bV^{\top} \bA = \UU^{\top}$, where $\bD^\dagger$ denotes the Moore--Penrose inverse of $\bD$. Consequently,
\[
\UU^{\top} x^+ = \bD^\dagger \bV^{\top} \bA x^+ = \bD^\dagger \bV^{\top} y,
\]
and hence $\bD^\dagger \bV^{\top} y$ is sparse as well. The Moore--Penrose inverse of $\bA$ is given by
$\bA^\dagger = \bar{\UU}\bD^{\dagger}\bV^{\top}.$ For noisy data $y^\delta = \bA x^+ + n$, the corresponding reconstruction is
\[
\bA^{\dagger}y^\delta = \bar{\UU}\bD^{\dagger}\bV^{\top} y^\delta = \bar{\UU}\left(\bD^\dagger \bV^{\top} y
+ \bD^{\dagger}\bV^{\top} n\right).
\]
This representation separates the reconstruction into the frame coefficients of the underlying signal and a noise term amplified by $\bD^{\dagger}$. Small entries of $\bD$ reflect the ill-conditioning of $\bA$ and can lead to strong amplification of the noise. The sparsity of the signal representation therefore provides a natural basis for regularization: rather than applying the inverse directly, the coefficients can be filtered to suppress noise while preserving the significant components of the signal. Such regularized reconstructions take the form
\begin{equation}\label{def:filter_reg}
    \mathcal{R}_\al (y^\delta) = \bar{\UU} f_\al\left(\bD^{\dagger},\bV^{\top} y^\delta\right),
\end{equation}
where the filter is applied component-wise. An example is the soft-thresholding filter \begin{equation}\label{def:soft_thresholding}
    f_\al(d,c) = \sign(c)\max\{|c|-\alpha/d,0\}.
\end{equation}
Under suitable assumptions on the filter and the parameter choice $\alpha=\alpha(\delta)$, such filtered reconstructions can yield a convergent regularization method in the sense of \cite{En96} as the noise level $\delta$ tends to zero, see \cite{ebner2023convergence, Fr19}.

Note that, although frame decompositions with non-tight frames are also used for inverse problems, non-tight frames in finite dimensions have the disadvantage that the frame bounds must be controlled and the corresponding dual frame must be computed.

\paragraph{Application 2: Sparse representation of observations.}
The matrix decomposition also provides a sparse representation of the observations. Indeed, from $\bA^{\top}\bV=\UU\bD$ we obtain $\bA=\bar{\bV}\bD\UU^{\top}$, and hence
\[
y=\bA x =\bar{\bV}\bD\UU^{\top} x.
\]
If the frame coefficients $\UU^{\top} x$ are sparse, then the coefficient vector
$\bD\UU^{\top} x=\bV^{\top} y$ is sparse as well. Thus, the observations $y$ admit a sparse representation with respect to the dual frame $\bar{\bV}$,
\[
y=\sum_{i=1}^K d_i \inner{x}{u_i}\bar v_i,
\]
where $\bar v_i$ denotes the $i$th element of the dual frame $\bar{\bV}$. Consequently, the observations can be represented and stored using only a small number of nonzero frame coefficients. In particular, instead of storing the full observation vector $y$, one may store its sparse frame representation $\bV^{\top} y=\bD\UU^{\top} x$ together with the corresponding indices.


\subsection{Main contributions}

In this paper, we develop a numerical methodology for solving the structured matrix decomposition (Problem~\ref{decom_problem}). Our approach combines an algebraic reduction of the decomposition problem with a new non-smooth and non-convex hierarchical (bilevel) optimization formulation of the original matrix decomposition problem. In particular, we develop from scratch a new bilevel optimization formulation on the Stiefel manifold containing the desired frame matrices. 

We first show that Problem~\ref{decom_problem} can be reduced to an equivalent problem involving only positive diagonal matrices for $\bA$. This reduction allows us to characterize necessary and sufficient conditions on the diagonal matrix $\bD$ for the existence of a decomposition satisfying the prescribed tightness and decomposition constraints. Once $\bD$ is fixed, the remaining problem can be formulated entirely in terms of the frame $\UU$, since the second frame $\bV$ is uniquely determined by $\UU$ and $\bD$. The tightness constraint on $\UU$ naturally leads to an optimization problem on the Stiefel manifold.

To incorporate the sparsity requirement, we formulate the resulting problem as a simple bilevel optimization problem on the Stiefel manifold \eqref{simple_bilevel}. The lower-level problem is a smooth but nonconvex optimization problem whose solution guarantees the tightness of $\bV$ and the decomposition property, while the upper-level objective is convex and lower semi-continuous, incorporating the sparsity requirement. The tightness of $\UU$ follows directly from the definition of the Stiefel manifold. Although the lower-level problem is convex when considered in the ambient Euclidean space, this convexity is generally lost after restricting the problem to the Stiefel manifold. We address the resulting nonconvex bilevel problem using a penalized subgradient method formulated in a general setting, which therefore applies beyond the specific matrix decomposition considered here.

We provide a rigorous convergence analysis of both levels of the resulting bilevel scheme. In particular, we establish convergence of the lower-level optimization procedure and analyze the convergence of the iterates to stationary points of the bilevel problem under the stated assumptions. The theoretical results are complemented by numerical experiments that demonstrate the practical behavior of the proposed method and its ability to identify sparse frame representations adapted to the given data.

\subsection{Related work}
Optimization-based approaches to frame computation have been studied, for example, by ~\cite{TroppIEEE05}, who develop convex optimization methods for constructing frames with prescribed structural properties. Our setting differs in that the frame is not designed independently of an operator; instead, the frame and a diagonal scaling are jointly required to represent a given matrix while promoting sparsity of its action on a data set.
A related line of work concerns scalable frames, where one seeks diagonal scalings that transform a given frame into a Parseval frame; see, e.g.,~\cite{KUTYNIOK20132225}.
Scalable-frame theory typically considers a fixed frame and asks whether a positive diagonal scaling yields a Parseval frame. In contrast, we optimize the frame itself together with the diagonal scaling, with the additional requirement that their product represent a prescribed operator and with a data-dependent sparsity objective.

Optimization subject to orthogonality constraints is naturally formulated on matrix manifolds. Structural optimization on embedded submanifolds was studied, for example, in~\cite{CambierAbsilSICOMP16}. For optimization on the Stiefel manifold, proximal-gradient methods have been developed in~\cite{ChenSIOPT2020}; see~\cite{Chen24} for a broader treatment of proximal and related methods for nonsmooth optimization on matrix manifolds. Weakly convex optimization, which is relevant to the nonsmooth structure arising in our analysis, has been studied extensively in~\cite{DavDrusSIOPT19,li2021weakly}.

Penalty methods are classical approaches to solve constrained optimization problems via a sequence of unconstrained problems. In convex optimization in Euclidean spaces, prox-penalization and related diagonal-gradient approaches have been studied extensively~\cite{Att96,BaiCom01,AttCzar18,AttCzarPey11,10530399,PeyStaQu26}. In contrast, penalty formulations specifically tailored to optimization on the Stiefel manifold are comparatively less developed. For example, a Lagrangian relaxation approach for optimization with Stiefel constraints is considered in~\cite{XiaPenaltySIOPT21}. Our work combines these ideas by employing a penalty formulation together with a subgradient method on the Stiefel manifold to solve the bilevel frame decomposition problem.


\section{Preliminaries}
\label{sec:Background}
%

Throughout, we work in real finite-dimensional vector spaces equipped with an inner product $\inner{\cdot}{\cdot}$ and the associated norm $\norm{\cdot}$. Unless stated otherwise, we use the Frobenius norm for matrices and the $\ell^2$-norm for vectors, i.e.,
\[
\norm{A}_F := \left(\sum_{i,j} |A_{ij}|^2\right)^{1/2},
\qquad
\norm{x}_2 := \left(\sum_i |x_i|^2\right)^{1/2},
\]
where we omit the subscript when it is clear from the context.

\subsection{Frames}\label{sec:frames}

We briefly recall some basic facts about finite-dimensional frame theory. More background on frame theory in Hilbert spaces can be found in \cite{Ch02}, while \cite{Waldron2018} provides an introduction to finite  tight frames. Let $R,K\in\mathbb{N}$ with $K\geq R$. A family $(u_i)_{i=1}^K\subset\mathbb{R}^R$ is called a \emph{frame} for $\mathbb{R}^R$ if there exist constants $A,B>0$ such that
\begin{equation}
A\norm{x}^2
\leq
\sum_{i=1}^K \abs{\inner{u_i}{x}}^2
\leq
B\norm{x}^2
\qquad\text{for all }x\in\mathbb{R}^R.
\end{equation}
In finite dimensions, we identify the frame with the matrix $\UU\in\mathbb{R}^{R\times K}$ whose columns are the frame elements $(u_i)_{i=1}^K$.
The constants $A$ and $B$ are called the lower and upper frame bounds of $\UU$, respectively. The frame is called \emph{tight} if $A=B$. A simple example of a tight frame is the union of $p$ orthonormal bases, in which case the frame bounds are $A=B=p$. Other prominent examples include finite-dimensional versions of wavelet, curvelet, and Gabor frames, which arise, for example, from suitable discretizations of these constructions in imaging applications \cite{candes2002recovering,candes2004new,daubechies1992ten}. A tight frame with $A=B=1$ is called a \emph{Parseval frame}.
Furthermore, a frame is called \emph{exact} if deleting any single frame element results in a family that is no longer a frame. A frame that is not a basis is called \emph{overcomplete}.

For a frame represented by $\UU\in\mathbb{R}^{R\times K}$, the corresponding analysis and synthesis operations are given by the matrix products
\[
\UU^\top x\in\mathbb{R}^K
\qquad\text{and}\qquad
\UU c\in\mathbb{R}^R,
\]
respectively. Thus, $\UU^\top$ maps a vector to its frame coefficients, while $\UU$ synthesizes a vector from a given coefficient vector. The associated frame operator is therefore
\[
\bS=\UU\UU^\top,
\]
which is symmetric positive definite, and hence invertible. In particular, the frame condition can equivalently be written as
$A\Id_{\mathbb{R}^R}\preceq \UU\UU^\top\preceq B\Id_{\mathbb{R}^R}$. For a tight frame with frame constant $A$, the frame condition is equivalent to $\UU\UU^\top=A\Id_{\mathbb{R}^R}$, while for a Parseval frame it is $\UU\UU^\top=\Id_{\mathbb{R}^R}$. Notably, every tight frame can be scaled to a Parseval frame: if $\UU$ is a tight frame with frame constant $A$, then $\UU/\sqrt{A}$ is a Parseval frame.

The redundancy of a frame generally leads to multiple possible coefficient representations of a given vector. In particular, there are different frames $\bar\UU\in\mathbb{R}^{R\times K}$ that can serve as dual frames and satisfy the reconstruction identity
\[
x=\bar\UU\UU^\top x
\qquad\text{for all }x\in\mathbb{R}^R,
\]
or, equivalently, $\bar\UU\UU^\top=\Id_{\mathbb{R}^R}$. 
A distinguished choice is the canonical dual frame, whose matrix is given by
$\bar\UU_{\mathrm{can}} = (\UU\UU^\top)^{-1}\UU$.
For a tight frame with frame constant $A$, this reduces to
$\bar\UU_{\mathrm{can}}=\frac{1}{A}\UU$ and for a Parseval frame, we simply have $\bar\UU_{\mathrm{can}}=\UU$.
Although the canonical dual is a natural choice, overcomplete frames generally admit multiple dual frames, and this freedom can be beneficial in applications.

\subsection{Stiefel manifold} \label{sec:Stiefel}

In the following sections, we formulate Problem~\ref{decom_problem} as an optimization problem over the Stiefel manifold. We therefore introduce the necessary notation here, while more details are collected the Appendix.

The Stiefel manifold for parameters $K,R$ is
\[
\St(K,R) \coloneqq  \{\mathbf{B} \in  \R^{K \times R} \mid \mathbf{B}^{\top} \mathbf{B} = \Id_{R} \}
\]
which is an embedded submanifold of the Euclidean space $\mathbb{R}^{K\times R}$. We use this convention for the Stiefel manifold, which is standard in the literature, naturally represents Parseval frames through their transpose: If $\UU\in\mathbb{R}^{R\times K}$ is a Parseval frame for $\mathbb{R}^R$, then $\UU^\top\in\St(K,R)$, and conversely, for every $\bB\in\St(K,R)$, the matrix $\bB^\top$ is a Parseval frame for $\mathbb{R}^R$.
The condition $\bB^\top\bB=\Id_R$ implies that the columns of $\bB$ are orthonormal, which yields
\[ \norm{\bB}_2=1,\qquad \norm{\bB}_F =\sqrt{R}, \qquad
\norm{b_{i}}_2 = 1,
\]
where $b_i$ are the columns of $\bB$.
Moreover, $\St(K,R)$ is compact, being a closed and bounded subset of $\mathbb{R}^{K\times R}$. If the dimensions are clear from the context, we simply write $\St$ for the Stiefel manifold.


\section{Computing data-aware tight frame decompositions}
\label{sec:Opt}
%

Let us recall the notion of diagonal frame decomposition (DFD) in the finite-dimensional setting. For a matrix $\bA \in \R^{N \times M}$, a DFD consists of a triplet $(\UU,\bV,\bD)$ such that, for some $K \geq \rank(\bA)$,
\begin{enumerate}
\item $\UU \in \R^{M \times K}$ is a frame for $\ker(\bA)^\perp$,
\item $\bV \in \R^{N \times K}$ is a frame for $\range(\bA)$,
\item $\bD = \diag(\dd_1,\ldots,\dd_K) > 0$ is diagonal and satisfies the quasi singular relation 
\[\bA^{\top}\bV = \UU\bD.\]
\end{enumerate}
We refer to a DFD as a \emph{tight DFD} if both $\UU$ and $\bV$ are tight frames, that is,
\[
\UU\UU^{\top} = A_{\UU}\Id_{\ker(\bA)^\perp}
\qquad\text{and}\qquad
\bV\bV^{\top} = A_{\bV}\Id_{\range(\bA)},
\UU
\]
where $A_{\UU}$ and $A_{\bV}$ denote the corresponding frame bounds.
Depending on the choice of $K$, the set of all tight DFDs can be large. To select a decomposition that is adapted to a given data set and promotes sparse representations, assume that we have collected data from the population of possible parameters $\bx$ stored in the data set
\[
\bX=[\bx_1,\ldots,\bx_S]\in\R^{M\times S}.
\]
To respect the problem structure captured in these sampled data points and to obtain an expressive dictionary for representing the ground truth parameter, a popular formulation of Problem~\ref{decom_problem} is the $\ell^1$-minimization problem
\begin{eqnarray}\label{eq:Opt}
\begin{aligned}
&\min_{\UU}\sum_{s=1}^{S}\norm{\UU^{\top}\bx_s}_{\ell^{1}(\R^{K})}\ \\[0.1cm]
&\text{subject to: } (\UU,\bV,\bD)\text{ is a tight DFD for $\bA$.}
\end{aligned}
\end{eqnarray}
Here, $\bA$ and the data set $\bX$ are given, while $\UU$, $\bV$, and $\bD$ are to be determined. Note that the objective function of the optimization problem is convex, but the constraints are not. Hence, in its current form, the optimization problem is challenging and it is not clear whether the above variational formulation can be exploited computationally for data-aware frame design. We address this question through a sequence of transformations that reduce the original problem to a more tractable formulation. In particular, we show that the problem can be recast as a nonconvex bilevel optimization problem over a Stiefel manifold. This reformulation enables the application of a first-order optimization method for which we establish convergence guarantees.

\subsection{Reduction to Parseval tight DFDs for positive diagonal matrices}

First, we note that the frames in a tight DFD can be rescaled arbitrarily, as the scaling can be absorbed into the quasi-singular values $\bD$. Indeed, if $(\UU,\bV,\bD)$ is a tight DFD of $\bA$, then
$(a \UU,b \bV,\frac{b}{a}\bD)$ is also a tight DFD of $\bA$ for any $a,b>0$. Since the $\ell^1$ sparsity term $\|\UU^\top \bx_s\|_{\ell^1}$ is not invariant under such rescaling, whereas the ideal $\ell^0$ sparsity measure is, we fix the scaling of the frames throughout. In particular, we w.l.o.g normalize the frames to be Parseval frames.

Our second reduction assumes that a reduced singular value decomposition (SVD) of the matrix $\bA\in\R^{N\times M}$ is available. Let $R=\rank(\bA)$ and consider $\bA=\bQ\Sigma\bW^{\top}$, 
where $\bQ\in\R^{N\times R}$ and $\bW\in\R^{M\times R}$ have orthonormal columns, and $\Sigma=\diag(\sigma_1,\ldots,\sigma_R)$ with $\sigma_1\geq\cdots\geq\sigma_R>0$.
Such a reduced SVD can be computed routinely using standard numerical linear algebra packages. The following lemma shows that, for a prescribed sequence of positive weights $(d_k)$, the construction of a DFD can be reduced to a problem involving only the positive diagonal matrix $\Sigma$.

\begin{lemma}
\label{lem:DFD_SVD_equivalence}
Let $\bA = \bQ \Sigma \bW^{\top}$ be the reduced SVD of $\bA \in \R^{N \times M}$, where $\bQ \in \R^{N \times R}$, $\bW \in \R^{M \times R}$ have orthonormal columns and $\Sigma = \diag(\sigma_1, \ldots, \sigma_R) \in \R_{+}^{R \times R}$ with $\sigma_1 \geq \cdots \geq \sigma_R > 0$ ($R=\rank(\bA)$). Let $\bD = \diag(d_1, \ldots, d_K)$ with $d_1 \geq \cdots \geq d_K > 0$ for some $K \geq R$.

\begin{enumerate}[label=(\roman*)]
\item \label{DFD_SVD:i} If $(\UU, \bV, \bD)$ is a DFD for $\bA$, then $(\tilde{\UU}, \tilde{\bV}, \bD)$, defined by
\begin{align*}
\tilde{\UU} &:= \bW^{\top} \UU, \\
\tilde{\bV} &:= \bQ^{\top} \bV,
\end{align*}
is a DFD for $\Sigma$. 
\item \label{DFD_SVD:ii} Conversely, if $(\tilde{\UU}, \tilde{\bV}, \bD)$ is a DFD for $\Sigma$, then $(\UU, \bV, \bD)$, defined by
\begin{align*}
\UU &:= \bW \tilde{\UU}, \\
\bV &:= \bQ \tilde{\bV},
\end{align*}
is a DFD for $\bA$. 
\end{enumerate}
In particular, $\UU$ and $\tilde{\UU}$, as well as $\bV$ and $\tilde{\bV}$, have the same frame bounds, respectively.
\end{lemma}

\begin{proof}
To prove \ref{DFD_SVD:i}, let $(\UU, \bV, \bD)$ be a DFD for $\bA$. Since $\bQ$ and $\bW$ have orthonormal columns spanning $\overline{\range(\bA)}$ and $\ker(\bA)^\perp$, respectively, $\tilde{\UU}$ and $\tilde{\bV}$ are again frames, and by Parseval's identity the frame bounds of $\UU$ and $\bV$ are preserved.

It remains to show that $(\tilde{\UU}, \tilde{\bV}, \bD)$ satisfies the quasi-singular relations for $\Sigma$. By the quasi-singular relation for $\bA$, we have $\bA^{\top} \bV = \UU \bD^{\top}$, and by the SVD, $\Sigma^{\top} \bQ^{\top} = \bW^{\top} \bA^{\top}$. Hence,
\[
\Sigma^{\top} \tilde{\bV}
= \Sigma^{\top} \bQ^{\top} \bV
= \bW^{\top} \bA^{\top} \bV
= \bW^{\top} \UU \bD^{\top}
= \tilde{\UU} \bD.
\]
Similarly, to prove \ref{DFD_SVD:ii} the same argument for the frame bounds holds and using $\bA^{\top} \bQ  = \bW \Sigma^{\top}$ and $\Sigma^{\top} \tilde \bV = \tilde \UU \bD^{\top}$, we obtain
\begin{equation*}
\bA^{\top} \bV = \bA^{\top} \bQ \tilde \bV
= \bW \Sigma^{\top} \tilde \bV
= \bW \tilde \UU\bD^{\top}
= \UU \bD^{\top}. \hfill \text{\qedhere}
\end{equation*}
\end{proof}

Lemma~\ref{lem:DFD_SVD_equivalence} allows us to reduce the construction of tight DFDs for a general matrix to the corresponding problem for a positive diagonal matrix $\Sigma\in\R^{R\times R}$, where the characterization of a tight DFD is particularly simple. Since $\Sigma$ has full rank, we have $\ker(\Sigma)^\perp=\operatorname{ran}(\Sigma)=\R^R$. Together with the normalization to Parseval frames, the frame conditions therefore reduce to $\tilde \UU\tilde \UU^{\top}=\Id_R$ and $\tilde \bV\tilde \bV^{\top}=\Id_R$. In summary, solving \eqref{eq:Opt} reduces to
\begin{eqnarray}\label{eq:Opt_parseval}
\begin{aligned}
&\min_{\tilde \UU}\sum_{s=1}^{S}\norm{\tilde \UU^{\top}\bW^{\top}\bx_s}_{\ell^{1}(\R^{K})}\ \\[0.1cm]
&\text{subject to: }  \tilde \UU \tilde \UU^{\top} = \tilde \bV \tilde \bV^\top = \Id_R \text{ and } \Sigma \tilde \bV = \tilde \UU \bD,
\end{aligned}
\end{eqnarray}
where $\Sigma$, $\bW$, $\bQ$, and $\bX$ are given, while $\tilde \UU$, $\tilde \bV$, and $\bD$ are to be determined. Then $(\bW\tilde \UU,\bQ\tilde \bV,\bD)$ yields the desired DFD for $\bA$.


\subsection{Predetermined quasi-singular values} \label{sec:Pre_D}

The next simplification is that we can eliminate the variable $\bV$, since any DFD  $(\UU,\bV, \bD)$ of $\Sigma$ satisfies $\bV=\Sigma^{-1} \UU\bD$. Hence, we can determine $\bV$ once we have computed $\tilde \UU$ (for a given $\bD$). Exploiting this linear coupling between the frames, we arrive at the simpler form
\begin{equation}\label{eq:min-U-s.t.-C}
        \min_{\UU}\sum_{s=1}^{S}\norm{ \UU^{\top} \bW^{\top}\bx_{s}}_{\ell^{1}(\R^{K})}\quad \text{s.t. } \UU\in\mathcal{C}_{\mathbf{D}} \coloneqq \left\{
\mathbf{U} \in \R^{R \times K} \;\middle|\;
\mathbf{U}\mathbf{U}^{\top} = \Id_R,
\quad
{\Sigma}^{-1}\mathbf{U}\mathbf{D}^2\mathbf{U}^{\top}{\Sigma}^{-1} = \Id_R
\right\}.
\end{equation}
If we fix the quasi-singular values $\bD$, the optimization reduces to a problem in the single variable $\UU$. It then remains to find a suitable strategy for choosing the quasi-singular values so that the objective in~\eqref{eq:min-U-s.t.-C} is as small as possible. This is a challenging problem for which no general purpose answer can be given. However, we can characterize the choices of $\bD$ for which a Parseval tight DFD exists, which gives us a numerical range of values to design potential candidates for the entries in the diagonal matrix $\bD$. The following proposition provides a necessary and sufficient feasibility condition and thus a potential criterion for selecting $\bD$. The result follows directly from the Fan--Pall theorem, stated in its original form for the real case in \cite[Theorem 1]{Fan_Pall_1957}.

\begin{theorem}[Fan--Pall] \label{thm:Fan_Pall}
Let $\bA \in \mathbb{R}^{K \times K}$ and $\bB \in \mathbb{R}^{R \times R}$ be two symmetric matrices with $K \geq R$. Let $\alpha_1 \geq \cdots \geq \alpha_K$ and $\beta_1 \geq \cdots \geq \beta_R$ denote the eigenvalues of $\bA$ and $\bB$, respectively. Then there exists a matrix $\bC \in \mathbb{R}^{K \times R}$ with $\bC^{\top} \bC =\Id_{R}$ such that
$\bC^{\top} \bA \bC = \bB$
if and only if
\[
\alpha_i \geq \beta_i \geq \alpha_{i+K-R}, \qquad i=1,\ldots,R.
\]
\end{theorem}

\begin{proposition}[Characteristic of the feasible set]\label{prop:feas_set}
Let $\Sigma = \diag(\sigma_1, \ldots, \sigma_R) \in \R^{R \times R}$ with $\sigma_1 \geq \cdots \geq \sigma_R > 0$ and $\bD = \diag(d_1, \ldots, d_K)$ with $d_1 \geq \cdots \geq d_K > 0$ for some $K \geq R$. 
A necessary and sufficient condition for the set $\scrC_\bD$ to be non-empty is the inequality
\begin{equation}\label{eq:feas_set}
    d_i \geq \sigma_i \geq d_{i+K-R} \quad \text{ for } i=1, \ldots, R.
\end{equation}
\end{proposition}
\begin{proof}
In our setting, $\bD^2$ and $\Sigma^2$ are not only symmetric but also diagonal and positive definite, with their eigenvalues given by the diagonal entries. Therefore, condition \eqref{eq:feas_set} is necessary and sufficient for the existence of a matrix $\UU \in \mathbb{R}^{R \times K}$ satisfying $\UU \UU^{\top} = \Id_R$ and $\UU \bD^2 \UU^{\top} = \Sigma^2$.
\end{proof}

Proposition~\ref{prop:feas_set} provides a necessary and sufficient condition for $\scrC_{\bD}\neq\emptyset$. In particular, if $K>R$, the construction of the scaling weights $(d_i)_{i=1}^K$ only requires checking $2R$ inequalities. Since $R$ is fixed for the problem under consideration, this verification has constant computational cost.

The choice of the number of scaling coefficients $K$ has a direct impact on the flexibility of the resulting frame decomposition. In the limiting case $K=R$, the Parseval constraint leaves essentially no redundancy in the frame $\UU$: A Parseval DFD only exists if $(d_i)_i=(\sigma_i)_i$, in which case $\UU$ can only be a signed permutation matrix and hence, can not be adapted to the data.

For $K>R$, the additional rows of $\UU$ provide additional degrees of freedom. This flexibility can be made precise by considering two choices $\tilde{\bD}\in\mathbb{R}^{\tilde K\times\tilde K}$ and $\bD\in\mathbb{R}^{K\times K}$, with $\tilde K<K$, such that the diagonal entries of $\tilde{\bD}$ form a subset of those of $\bD$. If $\tilde{\UU}\in\scrC_{\tilde{\bD}}$, then $\UU$ can be embedded into $\mathbb{R}^{K\times R}$ by adding $K-\tilde K$ zero rows (see Example~\ref{ex:DFD} below). The resulting matrix belongs to $\scrC_{\bD}$. Consequently, the feasible set for the larger value of $K$ contains an embedded copy of the feasible set for the smaller one. Thus, increasing $K$ can only improve the sparsity objective when additional scaling coefficients are added to $\bD$.
In practice, however, increasing $K$ beyond a certain point may provide little additional benefit, as the additional rows of $\UU$ may simply converge to zero rather than contributing meaningful additional structure to the representation.

\begin{example}\label{ex:DFD}
   We illustrate the flexibility of the frame decomposition.Let 
   $$
   \Sigma=\begin{pmatrix} 2 & 0 \\ 0 & 3 \end{pmatrix}\quad \textnormal{and}\quad
    \bD=\begin{pmatrix} 1 & 0 & 0 & 0 \\ 0 & 2 & 0 & 0 \\ 0 & 0 & 3 & 0 \\ 0 & 0 & 0 & 4 \end{pmatrix}. 
    $$
Then one can verify that the matrix 
 $$
\mathbf{U}_1=\begin{pmatrix} 1/2 & \sqrt{1/4+\sqrt{10}/8} & \sqrt{3/4 - 3\sqrt{10}/14} & \sqrt{-1/4 + 5\sqrt{10}/56} \\ \sqrt{\sqrt{10}/4-1/2} & -\sqrt{7/6- \sqrt{10}/3} & -\sqrt{1/14+\sqrt{10}/28} & \sqrt{11/42+\sqrt{10}/21} \end{pmatrix}
$$
defines a Parseval frame as well as
\[
\setlength{\arraycolsep}{3pt}
\bV_1=\Sigma^{-1}\mathbf{U}_1\bD
=\begin{pmatrix}
1/4 & \sqrt{1/4+\sqrt{10}/8} &
\sqrt{27/16-27\sqrt{10}/56} &
\sqrt{-1+5\sqrt{10}/14}\\
\sqrt{\sqrt{10}/36-1/18} &
-\sqrt{14/27-4\sqrt{10}/27} &
-\sqrt{1/14+\sqrt{10}/28} &
\sqrt{88/189+16\sqrt{10}/189}
\end{pmatrix}.
\]
On the other hand, consider 
$$
\mathbf{U}_2=\begin{pmatrix} 0 & 1 & 0 & 0 \\ 0 & 0 & 1 & 0 \end{pmatrix} \text{ and } \mathbf{V}_2 =\Sigma^{-1}\mathbf{U}_2\bD=\begin{pmatrix} 0 & 1 & 0 & 0 \\ 0 & 0 & 1 & 0 \end{pmatrix}.
$$
These also define a Parseval DFD.
Thus, for the same $\Sigma$ and $\bD$, there are multiple Parseval DFDs. In particular, the second choice selects the scaling coefficients $2$ and $3$, which exactly match the singular values of $\Sigma$, whereas the first choice distributes the contributions across all four scaling coefficients. This illustrates the additional flexibility provided by $K>R$.
\end{example}

The preceding observations show that the feasibility condition alone does not determine a suitable choice of $\bD$. In particular, while $\bD$ must capture the spectral variation of $\Sigma$ sufficiently well to ensure a small value of the lower-level objective, it should at the same time provide sufficient freedom for $\UU$ to adapt to the data. This leads to a trade-off between spectral compatibility and the sparsity-promoting objective in~\eqref{eq:min-U-s.t.-C}. Moreover, the constrained problem over $\scrC_{\bD}$ remains potentially infeasible and computationally challenging. In the next subsection, we therefore relax the exact Parseval constraint for $\bV$. This relaxation makes every choice of $\bD$ admissible and, at the same time, allows for a more flexible optimization over $\UU$. We then introduce a data-driven heuristic for selecting the diagonal entries of $\bD$.

\subsection{Bilevel Optimization on the Stiefel Manifold}

 To obtain a numerically viable solution strategy, we adopt a relaxation approach in which the potentially inconsistent constraint $\mathbf{U}\in\mathcal{C}_{\mathbf{D}}$ is replaced by the best approximation 
\[
\left\{\UU \, \Big| \, \UU^{\top} \in \argmin_{\mathbf{U}\in\St(K,R)}\norm{{\Sigma}^{-1}\mathbf{U}^{\top}\mathbf{D}^2\mathbf{U}{\Sigma}^{-1} - \Id_R}_F^2 \right\}.
\]

Recall from Section~\ref{sec:Stiefel} that the set of all transposed Parseval frames can be identified with the Stiefel manifold $\St(K,R)$. In the following, we use $\UU\in\St(K,R)$ as the optimization variable, with the corresponding Parseval frame given by $\UU^{\top}$. With this notational convention, we formulate the problem using an outer approximation of the feasibility constraint. This leads to the following \emph{simple bilevel optimization problem}~\citep{dempe2021simple} over the Stiefel manifold:
\begin{equation} \label{pen_opt}
    \min_{\UU\in\mathcal{S}} \;\phi(\UU)
    \coloneqq \|\UU \bW^{\top}\bX\|_1,
    \qquad
    \mathcal{S}\coloneqq\argmin_{\UU\in\St(K,R)}\Psi(\UU) \coloneqq\left\| \Sigma^{-1}\UU^{\top}\bD^2\UU\Sigma^{-1}-\Id_R\right\|_F^2.
\end{equation}
The upper level problem is the data-driven sparsity term in order to inject into the learned frames  the structure of the data underlying. Our bilevel formulation of the frame learning problem allows us to perform our optimization even if the scaling factors $\bD$ are badly chosen, so that $\scrC_{\bD}=\emptyset$. Indeed, if $\bD$ is chosen so that $\scrC_{\bD} \neq \emptyset$, then $\Psi(\scrS)=0$ and the original problem is solved exactly. If $\scrC_{\bD} \neq \emptyset$, we still obtain a DFD in which $\UU$ is a Parseval frame, and the frame bounds of $\mathbf{V}$ can be controlled in terms of the value of the lower-level objective function on the solution set. 

\begin{proposition}
Let $\UU \in \St(K,R)$ and define $\bV = {\Sigma}^{-1}\UU^{\top} \bD$. Suppose $\Psi(\UU) = \epsilon < 1$. Then $(\UU^{\top}, \bV, \bD)$ is a DFD for $\Sigma$, where $\UU^{\top}$ is a Parseval frame and $V$ has lower and upper frame bounds $1- \epsilon$ and $1+\epsilon$, respectively.
\end{proposition}

\begin{proof}
$\UU^{\top}$ is a Parseval tight frame. By the fact that the spectral norm is bounded above by the Frobenius norm and the definition of $\mathbf{V}$, we have 
\[\norm{\Id - \bV \bV^{\top}}_2 \leq \norm{\Id - \bV \bV^{\top}}_F = \Psi(\UU) = \epsilon
\]
hence, the eigenvalues of $\bV \bV^{\top}$ are in the interval $[1-\epsilon, 1+ \epsilon]$. 
\end{proof}

For the application to filtered reconstruction of inverse problems, it is not essential that $\bV$ be a Parseval or tight frame. Although computing its canonical dual would be more involved, this is not required in our reconstruction method, where only the dual of $\UU$ is needed. Nevertheless, the stability of the reconstruction depends, among other factors, on the frame bounds, or equivalently the condition number, of $\bV$. These need not be equal to one, but should remain reasonably well-conditioned.
A key requirement for a DFD is that the ill-conditioning of $\Sigma$, and hence of $\bA$, is captured by $\bD$, since the corresponding components are subsequently filtered. To illustrate this, consider the extreme choice $\bD=\Id_K$. For any Parseval frame $\UU$, the corresponding frame $\bV = \Sigma^{-1} \UU$ inherits the ill-conditioning of $\Sigma$; in particular, its condition number can be as large as that of $\Sigma$. Thus, the primary objective in choosing $\bD$ is to ensure that $\min\Psi$ is zero or sufficiently small, which generally requires $\bD$ to capture the dominant spectral variation of $\Sigma$. Among such choices, $\bD$ should then be selected to promote sparsity in $\UU$.

A heuristic way to compute the diagonal entries of $\bD$ that incorporates the structure of the data is as follows. Let $\bX\bX^\top = Q_X \Lambda Q_X^\top$ be the eigendecomposition of $\bX \bX^\top$. We define
\[
\hat{\dd} = (\hat d_i)_{i=1}^K = \argmin_{\dd} \left\| \Sigma^{-1} \bW^{\top} Q_X \diag(\dd)^2 Q_X^\top \bW \Sigma^{-1} - \Id_R \right\|_F^2 
\] 
and choose $\bD = \diag(d_1, \ldots, d_K)$ with 
\begin{equation}\label{eq:d}
d_i = \begin{cases}
 \quad \hat d_i &\text{ where } \hat d_i> \epsilon \\[0.1cm]
\norm{\Sigma}_F^2/R \quad &\text{ otherwise.}
\end{cases}
\end{equation}
Thus, the resulting scaling coefficients are adapted to the dominant directions of the data while simultaneously taking the spectral scaling induced by $\Sigma$ into account. In particular, this choice aims to align the frame \(\UU\) with directions that are relevant to the data while keeping the lower-level objective \(\Psi\) small, but not necessarily zero.
In practice, some of the coefficients $\hat d_i$ may become very small or vanish. Such coefficients provide little additional freedom for the frame and may lead to degenerate frame elements. We therefore replace coefficients below a prescribed threshold $\epsilon$ with a neutral scaling value.

Although \eqref{eq:d} favors sparsity in the frame \(\UU\), it may still lead to large values of the lower-level objective \(\Psi\) and, consequently, a high condition number of the resulting frame \(\bV\). This can potentially cause instability when the frame is used for regularization in inverse problems. As an alternative, the diagonal entries can be chosen independently of the data using the following construction, which is designed to satisfy the conditions of the Fan--Pall theorem. The practical implications of this choice depend on the application and will be investigated numerically in Section~\ref{sec:Num}.
\begin{equation}
d_i^2 =
\begin{cases}
\displaystyle
\frac{1}{i}\sum_{j=1}^{i}\sigma_j^2,
& 1\leq i\leq K-R+1 \wedge i \leq R,
\\[2ex]
\displaystyle
\frac{1}{K-R+1}
\sum_{j=i-K+R}^{i}\sigma_j^2,
& K-R+2\leq i\leq R,
\\[2ex]
\displaystyle
\frac{1}{K-i+1}
\sum_{j=i-K+R}^{R}\sigma_j^2,
& R+1\leq i\leq K.
\end{cases}
\label{eq:D_fan_pall_short}
\end{equation}


\section{The Numerical Scheme}
\label{sec:algorithm}
%
In this section, we describe a general numerical scheme for computing data-aware DFDs. Since our complexity statements apply to a rather broad class of optimal design problems, we keep the presentation at an abstract level. Throughout this section, we use both Euclidean and Riemannian derivatives. To distinguish between the two, we denote Riemannian quantities by the subscript $\St$. The relevant tools for optimization on the Stiefel manifold are summarized in Appendix \ref{app:Riemann}.

Let $\mathcal{U} \subset \R^{K \times R}$ be an open, convex, and bounded neighborhood of $\St(K,R)$. We consider the simple bilevel optimization problem
\begin{equation} \label{simple_bilevel}
\min_{\UU \in \mathcal S} \phi(\UU),
\qquad 
\mathcal S \coloneqq \argmin_{\UU \in \St(K,R)} \Psi(\UU),
\end{equation}
where the functions $\phi : \mathcal U \to \R$ and $\Psi : \mathcal U \to \R$ satisfy the following assumptions:
\begin{enumerate}
    \item $\phi$ is convex and lower semicontinuous;
    \item $\Psi$ is continuously differentiable with a Lipschitz continuous Euclidean gradient:
    \[
    \norm{\nabla\Psi(\UU_{1})-\nabla\Psi(\UU_{2})}
    \leq L_{\nabla \Psi} \norm{\UU_{1}-\UU_{2}}
    \qquad \text{for all } \UU_{1},\UU_{2}\in\mathcal U.
    \]
\end{enumerate}
The lower-level problem determines the feasible set $\mathcal S$ through the minimizers of the smooth function $\Psi$ on the Stiefel manifold, while the upper-level objective $\phi$ selects a preferred solution among these minimizers. This model formulation allows us to approach a solution of Problem \ref{decom_problem} by identifying the data as in \eqref{pen_opt}. 

Note that, since the Stiefel manifold is not convex, the restriction of $\phi$ to $\St$ is generally not convex. However, by compactness of $\St$ and convexity of $\phi$ on $\mathcal U$, $\phi|_{\St}$ is Lipschitz continuous with some constant $L_\phi$. The same holds for the restriction $\Psi|_{\St}$, which is Lipschitz continuous with some constant $L_\Psi$. Moreover, since $\Psi$ has a $L_{\nabla\Psi}$-Lipschitz continuous Euclidean gradient on $\mathcal U$, it is $L_{\nabla\Psi}$-weakly convex on $\mathcal U$. We recall that a function $f\colon\scrU\to\R$ is called $\tau$-weakly convex if $f+\tau\norm{\cdot}^2/2 $ is convex for some $\tau\geq0$.
Furthermore, since the Stiefel manifold $\St$ is compact and $\phi$ is convex on the open set $\mathcal U$, the restriction $\phi|_{\St}$ is continuous and hence attains both its minimum and maximum on $\St$. Likewise, since $\Psi$ is continuously differentiable on $\mathcal U$, both $\Psi|_{\St}$ and $\|\nabla\Psi\|_{\St}$ attain their respective minimum and maximum on $\St$. In particular, also the Riemannian gradient $\nabla_{\St}\Psi$ is Lipschitz continuous on $\St$ with the constant
\[L_{\nabla_{\St}\Psi}=L_{\nabla\Psi}+2\max_{\UU\in\St}\|\nabla\Psi(\UU)\|_F.
\]

\begin{remark}
    Note that the $\Psi$ defined in \eqref{pen_opt} can be extended to the Euclidean space $\R^{K\times R}$ but fails to be Lipschitz smooth on the full domain. However, the claimed Lipschitz smoothness holds on any bounded neighborhood of the Stiefel manifold $\St$. 
\end{remark}

\subsection{Penalty Subgradient Method on the Stiefel Manifold}
%
To approximate a solution of problem \eqref{simple_bilevel}, we extend recent algorithmic advances for weakly convex minimization over the Stiefel manifold \cite{li2021weakly} to simultaneously address both levels of the problem. To this end, we adopt the generalized Tikhonov regularization approach of \cite{solodov2007explicit} and consider the time-varying auxiliary objective function
\begin{equation}
    \Omega_{k}(\UU):=\phi(\UU)+\beta_{k}\Psi(\UU):\scrU\to\R
\end{equation}
where $(\beta_{k})_{k}$ is non-decreasing penalty sequence, formally going to infinity.
By Theorem 4.1 in \cite{yang2014optimality}, this is a regular function, and thus  $\partial_{\St}\Omega_{k}(\UU)=\partial_{\St}\phi(\UU)+\beta_{k}\nabla_{\St}\Psi(\UU)$. 
The penalty term increasingly promotes minimization of the lower-level objective $\Psi$, while the term $\phi$ selects a preferred solution among the lower-level minimizers. 
This auxiliary function is injected into a Riemannian subgradient-type method, analyzed in \cite{li2021weakly} for the stationary case ($\beta_k \equiv \beta$). A single iteration reads as 
\begin{equation}\label{iterates}
\UU^{k+1}=\proj_{\St}\Bigl(\UU^k-\ga_k\bigl(g^k+\beta_k\nabla_{\St}\Psi(\UU^k)\bigr)
\Bigr),\qquad g^k \in \partial_{\St}\phi(\UU^k),
\end{equation}
where $(\ga_k)_k$ denotes the sequence of step sizes. The Riemannian gradient and subgradient ensure that the search direction lies in the tangent space of the Stiefel manifold, while the orthogonal projection maps the updated point back onto the manifold and thus serves as a retraction. This procedure is briefly discussed below and more detailed in the Appendix. The precise formulation of the resulting method is given in Algorithm~\ref{alg:RPS}. The convergence properties of the method are established in the subsequent section.

\begin{algorithm}[t]
\caption{Riemannian Penalty Subgradient Method}
\label{alg:RPS}
\begin{algorithmic}[1]
\Require $\UU^0 \in \St(K,R)$, step sizes $(\ga_k)_k$, penalty parameters $(\beta_k)_k$
\For{$k=0,1,\dots$}
    \State Evaluate the Euclidean gradient $\nabla \Psi(\UU^k)$ and select $\tilde g^k \in \partial \phi(\UU^k)$.
    \State Compute the Riemannian gradient and subgradient: 
    \[\nabla_{\St}\Psi(\UU^k)=\proj_{T_{\UU^k}\St}\bigl(\nabla \Psi(\UU^k)\bigr)\] 
    \[g^k=\proj_{T_{\UU^k}\St}(\tilde g^k)\in\partial_{\St}\phi(\UU^k) 
    \hfill\displaystyle\proj_{T_{\bA}\St}(\bV)=\bV-\bA\frac{\bA^\top\bV+\bV^\top\bA}{2}\]
    \State Take a tangent-space step: 
    \[\UU_{+}^{k}=\UU^k-\ga_k\bigl(g^k+\beta_k\nabla_{\St}\Psi(\UU^k)\bigr)
    \]
    \State Retract onto the Stiefel manifold:
    \[\UU^{k+1}=\proj_{\St}(\UU^{k}_+) \hfill\displaystyle\proj_{\St}(\bV)=\bQ\bW^{\top} \text{ with SVD } \bV=\bQ\Lambda\bW^{\top} 
    \]
\EndFor
\end{algorithmic}
\end{algorithm} 

\paragraph{Gradients and subgradients on the Stiefel manifold.}
Since $\St(K,R)$ is an embedded submanifold of $\R^{K\times R}$, the Riemannian gradient of $\Psi$ can be obtained by projecting its Euclidean gradient onto the tangent space
\[
\nabla_{\St}\Psi(\UU)=\proj_{T_{\UU}\St}\bigl(\nabla\Psi(\UU)\bigr).
\]
For the generally nondifferentiable convex function $\phi$, the projection of its Euclidean subdifferential onto the tangent space gives the Clarke subdifferential of its restriction to the Stiefel manifold
\[
\partial_{\St}\phi(\UU)=\proj_{T_{\UU}\St}\bigl(\partial\phi(\UU)\bigr).
\]
This follows from \cite[Thm.~5.1]{yang2014optimality}, since the restriction of $\phi$ to the Stiefel manifold is Lipschitz continuous.

\paragraph{Retractions on the Stiefel manifold.}
Since tangent vectors do not generally remain on the manifold after an explicit update step, a retraction is used to map the resulting point back onto $\St$. In this work, we employ the polar decomposition-based retraction to exploit its nonexpansiveness.

\subsection{Convergence Analysis}

Since the bilevel design problem \eqref{simple_bilevel} is generally nonsmooth and nonconvex on the Stiefel manifold, convergence to global minimizers cannot be expected in general. Instead, we establish asymptotic convergence to suitably defined stationary points. To quantify stationarity, we follow the approach of \cite{li2021weakly} and use a Moreau-type envelope. For $\lambda>0$ and a $\tau$-weakly convex function
$h:\scrU\to\mathbb{R}$, we define
\begin{equation}
\label{eq:moreau}
\left\{
\begin{aligned}
h_\lambda(\mathbf{U})& = \min_{\mathbf{V}\in \St(K,R)}\left\{
h(\mathbf{V})+\frac{1}{2\lambda}\|\mathbf{V}-\mathbf{U}\|_F^2\right\},
\qquad\mathbf{U}\in \St(K,R),\\[0.3cm]
\prox_{\la h}(\mathbf{U})& = \operatorname*{argmin}_{\mathbf{V}\in \St(K,R)}\left\{
h(\mathbf{V})+\frac{1}{2\lambda}\|\mathbf{V}-\mathbf{U}\|_F^2\right\},
\qquad\mathbf{U}\in \St(K,R).
\end{aligned}
\right.
\end{equation}
Note that the proximal mapping need not be single-valued in general. 
However, the particular choice of $\hat \UU \in \prox_{\lambda h}(\UU)$ does
not affect the subsequent analysis, since every proximal point
satisfies the same first-order optimality condition. Throughout, whenever $\prox_{\lambda h}(\UU)$ is used as a point, we refer to an arbitrary selection from the proximal mapping.
For $\hat\UU \in \prox_{\lambda h}(\UU)$, the first-order optimality condition gives
\begin{equation*}
  \label{eq:first-order-opt-prox}
  0 \in \partial_{\mathcal \St}h(\hat \UU)+\frac{1}{\lambda}
\proj_{T_{\hat \UU}\St}
\bigl(\hat \UU-\mathbf{U}\bigr).
\end{equation*}
Consequently,
\begin{equation}
\label{eq:theta}
\dist(0, \partial_{\St} h(\hat \UU))
\leq\lambda^{-1}\left\|\proj_{T_{\hat \UU}\St}
\bigl(\hat \UU-\mathbf{U}\bigr)\right\|_F
\leq\lambda^{-1}\|\hat \UU-\mathbf{U}\|_F
=:\Theta_{\la h}(\mathbf{U}).
\end{equation}

Thus, $\Theta_{\lambda h}$ provides a stationarity measure for $h$. In particular,
$\Theta_{\lambda h}(\UU)=0$ implies $\hat\UU=\UU$ and hence
$0\in\partial_{\St}h(\UU)$.
This framework is applied to both levels of the bilevel problem. For the lower-level objective, we use $\Theta_{\lambda\Psi}$,
which measures violation of lower-level stationarity. 
In particular, since $\Psi$ is differentiable with Lipschitz continuous Riemannian gradient, we have
\begin{equation}\label{eq:riemann}
    \norm{\nabla_{\St}\Psi(\UU)}_F \leq L_{\nabla_{\St} \Psi} \snorm{\UU - \hat\UU}_F + \norm{\nabla_{\St}\Psi(\hat\UU)}_F
    \leq (L_{\nabla_{\St}\Psi} \la + 1) \, \Theta_{\la \Psi}(\UU).
\end{equation}
For the penalized objective $\Omega_k=\phi+\beta_k\Psi$, we consider
$\beta_k\Theta_{\lambda\Omega_k/\beta_k}$, as a measure of upper-level stationarity. 
In particular, \eqref{eq:theta} yields the bound
\[
\dist\bigl(0,\partial_{\St}\Omega_k(\prox_{\lambda\Omega_k/\beta_k}(\UU))\bigr)
\leq
\beta_k\Theta_{\lambda\Omega_k/\beta_k}(\UU).
\]

Another key technical tool in our analysis is the following Riemannian subgradient inequality on the Stiefel manifold for weakly convex functions, taken from \cite[][Theorem 1]{li2021weakly}. 

\begin{theorem}[Riemannian sugradient inequality]\label{thm:subgrad_ineq}
Let $\scrU \subset \R^{K \times R}$  be any bounded, open, and convex set that contains $\St(K,R)$. 
Suppose that $h \colon \scrU \rightarrow \mathbb{R}$ is $\tau$-weakly convex for some $\tau \geq 0$. Then, there exists a constant $L>0$ such that $h$ is $L$-Lipschitz continuous on $\St(K,R)$ and satisfies
\[h(\mathbf{\bV})\geq h(\mathbf{\UU})+\Bigl\langle\widetilde{\nabla}_{\mathcal \St} h(\mathbf{U}),\mathbf{V}-\mathbf{U}\Bigr\rangle-\frac{\tau + L}{2}\|\mathbf{V}-\mathbf{U}\|_F^2\]
for all
\[\widetilde{\nabla}_{\mathcal \St} h(\mathbf{U})\in\partial_{\mathcal \St} h(\mathbf{U}),\qquad\mathbf{U},\mathbf{V}\in \St(K,R).
\]
\end{theorem}

Compared to the corresponding Euclidean inequality for weakly convex functions \eqref{eq:weakly}, the Riemannian inequality contains the additional term $L$, where $L$ is a Lipschitz constant of $h$ on $\St$. This additional term arises from restricting $h$ to the nonconvex Stiefel manifold and accounts for the geometry of the manifold. For $\tau=0$, the inequality reduces to the corresponding Riemannian subgradient inequality for convex functions. In the smooth case, the Riemannian subgradient is replaced by the Riemannian gradient. Thus, the inequality provides a unified first-order estimate for the smooth lower-level objective and the potentially nonsmooth penalized objectives considered here.

We first investigate the decrease of the stationarity measure associated with the lower-level problem. Recall that  $\phi|_{\St}$ is $L_\phi$-Lipschitz continuous, $\Psi|_{\St}$ is $L_\Psi$-Lipschitz continuous, and $\Psi$ is $L_{\nabla\Psi}$-weakly convex on $\mathcal U$.

\begin{proposition}\label{prop:Theta_esimate}
Let $\{\mathbf U^k\}$ be the sequence generated by \eqref{iterates} with arbitrary initialization, step sizes $(\ga_k)_k$ and penalty parameters $(\beta_k)_k$.
Then, for every $\lambda < (L_{\nabla \Psi} + L_\Psi + L_\phi/\beta_0)^{-1}$,
the following estimate holds for all $k\geq0$:
\begin{equation}
\label{eq:prop-descent}
\la^2\Theta^2_{\la \Psi}(\mathbf U^k)
\leq \left(\frac{b_k + \sqrt{b_k^2- 4 a_k c_k}}{2 \abs{a _k}}\right)^2 \leq \frac{b_k^2}{a_k^2} + \frac{2c_k}{\abs{a_k}}
\end{equation} 
where
\begin{align}
    a_k &= \frac{\ga_k}{2\la}\left(L_\phi + \beta_k(L_{\nabla \Psi} + L_\Psi - \frac{1}{\la})\right) \label{eq:a}\\
    b_k &= \frac{\ga_k}{\la} L_\phi\label{eq:b}\\
    c_k &= (\Psi_\la(\UU^k) - \Psi_\la(\UU^{k+1})) + \frac{\ga_k^2}{2 \la}(\beta_k L_{\Psi} + L_\phi)^2.\label{eq:c}
\end{align}
\end{proposition}

\begin{proof}
By definition of the Moreau envelope, the operator $\prox_{\la \Psi}(\mathbf U^{k+1})$ and the fact that $\prox_{\la \Psi}(\mathbf U^k)\in\St(K,R)$, we have 
\begin{align}\label{eq:moreau-step}
\Psi_\lambda(\mathbf U^{k+1})
&=\Psi(\prox_{\la \Psi}(\mathbf U^{k+1}))+\frac{1}{2\lambda}
\|\prox_{\la \Psi}(\mathbf U^{k+1})-\mathbf U^{k+1}\|_F^2\\
&\leq 
\Psi(\prox_{\la \Psi}(\mathbf U^k))+\frac{1}{2\lambda}
\|\prox_{\la \Psi}(\mathbf U^k)-\mathbf U^{k+1}\|_F^2 .
\end{align}
The non-expansiveness condition eq. \eqref{eq:Retraction} gives
\begin{align*}
\|\prox_{\la \Psi}(\mathbf U^k)-\mathbf U^{k+1}\|_F^2
&\leq\|\mathbf{U}^k - \ga_k \big( g^k + \beta_k \nabla_{\St} \Psi(\mathbf{U}^k)\big)-\prox_{\la \Psi}(\mathbf U^k)\|_F^2\\
&=\|\mathbf U^k-\prox_{\la \Psi}(\mathbf U^k)\|_F^2 - 2 \ga_k \inner{\mathbf{U}^k-\prox_{\la \Psi}(\mathbf U^k)}{g^k + \beta_k \nabla_{\St} \Psi(\mathbf{U}^k)} \\
& \quad + \gamma_k^2\norm{g^k + \beta_k \nabla_{\St} \Psi(\mathbf{U}^k)}^2_{F}.
\end{align*}
To simplify the notation, let us set $\hat{\UU}^{k}\equiv \prox_{\lambda\Psi}(\UU^{k})$. Since the mapping $\Omega_k = \phi+\beta_{k}\Psi$ is $\beta_k L_{\nabla \Psi}$-weakly convex on $\mathcal{U}$ and $L_\phi + \beta_k L_\Psi$-Lipschitz continuous on $\St$, it follows 
\[
\inner{g_{k}+\beta_{k}\nabla_{\St}\Psi(\UU^{k})}{\UU^{k}-\hat{\UU}^{k}}\geq\Omega_{k}(\UU^{k})-\Omega_{k}(\hat{\UU}^{k})-\frac{\beta_{k}(L_{\nabla \Psi}+L_{\Psi})+L_{\phi}}{2}\norm{\UU^{k}-\hat{\UU}^{k}}^{2}_{F}. 
\]
Moreover, the Lipschitz continuity of $\phi$ and $\Psi$ implies
\begin{align}
\norm{g^k + \beta_k \nabla_{\St} \Psi(\mathbf{U}^k)}_{F}&\leq \beta_k L_\Psi + L_\phi,\text{ and }\\ 
\phi(\UU^{k})-\phi(\hat{\UU}^{k})&\geq -L_{\phi}\norm{\UU^{k}-\hat{\UU}^{k}}_{F}.
\end{align}
Using these bounds, we can continue the estimation with 
\begin{align*}
\|\hat{\UU}^{k}-\mathbf U^{k+1}\|_F^2&\leq\|\mathbf U^k-\hat{\UU}^{k}\|_F^2-2\gamma_k\bigl(\Omega_{k}(\mathbf U^k)-\Omega_{k}(\hat{\UU}^{k})\bigr)\\
&\quad+\gamma_k(\beta_k(L_{\nabla \Psi} + L_\Psi) + L_\phi)\|\mathbf U^k-\hat{\UU}^{k}\|_F^2+\gamma_k^2(\beta_k L_\Psi + L_\phi)^2\\
& \leq \left( 1 + \gamma_k(\beta_k(L_{\nabla \Psi}+ L_\Psi) + L_\phi)\right) \|\mathbf U^k-\hat{\UU}^{k}\|_F^2 \\
& \quad + 2\ga_k L_\phi \|\mathbf U^k-\hat{\UU^k}\|_{F}\\
&\quad  + 2 \ga_k \beta_k(\Psi(\hat{\UU}^k)-\Psi(\mathbf U^k))+ \gamma_k^2(\beta_k L_\Psi + L_\phi)^2
\label{eq:almost-final}
\end{align*}
Using the definition of the Moreau envelope once more, we see
\[
\Psi_\lambda(\mathbf U^k)=\Psi(\hat{\UU}^k)
+\frac{1}{2\lambda}
\|\mathbf U^k-\hat{\UU}^k\|_F^2 ,
\]
and since $\Psi_\lambda(\mathbf U^k)\leq \Psi(\mathbf U^k)$, we infer
\[
2\gamma_{k}\beta_{k}(\Psi(\hat{\UU}^{k})-\Psi(\UU^{k}))\leq-\frac{\gamma_{k}\beta_{k}}{\lambda}\norm{\UU^{k}-\hat{\UU}^{k}}^{2}_{F}.
\]
Injecting this into the penultimate bound, we are left with 
\begin{align*}
    \|\hat{\UU}^{k}-\mathbf U^{k+1}\|_F^2&\leq\left( 1 + \gamma_k(\beta_k(L_{\nabla \Psi}+ L_\Psi) + L_\phi)-\frac{\gamma_{k}\beta_{k}}{\lambda}\right) \|\mathbf U^k-\hat{\UU}^{k}\|_F^2 \\
& \quad + 2\ga_k L_\phi \|\mathbf U^k-\hat{\UU^k}\|_{F} +\gamma_k^2(\beta_k L_\Psi + L_\phi)^2
\end{align*}

Combining all these estimates,
\begin{align*}
  \Psi_\lambda(\mathbf U^{k+1})&\leq \Psi_\lambda(\mathbf U^{k})+\frac{\gamma_{k}}{2\lambda}\left[L_{\phi}+\beta_{k}\left(L_{\nabla \Psi}+L_{\Psi}-\frac{1}{\lambda}\right)\right]\norm{\UU^{k}-\hat{\UU}^{k}}^{2}_{F}\\
  &+ \frac{\gamma_{k} L_{\phi}}{\lambda}\norm{\UU^{k}-\hat{\UU}^{k}}_{F}+\frac{\gamma_{k}^{2}}{2\lambda}(\beta_{k}L_{\Psi}+L_{\phi})^{2} \\
 & = a_k \|\mathbf U^k-\hat{\UU}^{k}\|_F^2 + b_k \|\mathbf U^k-\hat{\UU}^{k}\|_{F} + c_k + \Psi_\lambda(\mathbf U^{k+1}),
  \end{align*}
where $a_{k},b_{k}$ and $c_{k}$ are defined in equations \eqref{eq:a}, \eqref{eq:b} and \eqref{eq:c}, respectively.
By the assumption on $\lambda$, we have $a_k<0$, which yields the desired result.  
The second estimate follows from the concavity of the square root function.
\end{proof}

\begin{theorem}\label{conv_stationary}
Let $\{\mathbf U^k\}$ be the sequence generated by \eqref{iterates} with arbitrary initialization, step sizes $(\ga_k)_k$ and monotonically increasing penalty parameters $(\beta_k)_k$. Then, for every $\lambda < (L_{\nabla \Psi} + L_\Psi + L_\phi/\beta_0)^{-1}$, the following estimate holds
\begin{equation}\label{eq:conv_stationary}
    \underset{0\leq k\leq T}{\min} \Theta^2_{\la \Psi}(\mathbf U^k) \leq \frac{4 \la L_\la \Delta_{\Psi}^{0} + 2 L_\la ( L_\Psi+ L_\phi/\beta_0)^2 \sum_{k=0}^T (\ga_k \beta_k)^2 + 4 L_\phi^2 \sum_{k=0}^T \ga_k/\beta_k}{\la^2 L_\la^2 \sum_{k=0}^T \ga_k \beta_k},
\end{equation}
where $L_\la \coloneqq 1/\la - L_{\nabla \Psi} - L_\Psi - L_\phi/\beta_0>0$ and $\Delta_{\Psi}^{0} \coloneqq \Psi(\UU^{0})-\min_{\St}\Psi$. 
In particular, the right-hand side of \eqref{eq:conv_stationary} converges to zero as $T \to \infty$ if and only if the following design criteria for the step-size and penalty sequences are satisfied:
\[
\sum_{k=0}^{T}\gamma_k\beta_k \to\infty, \qquad
\frac{\sum_{k=0}^{T}(\gamma_k\beta_k)^2}{\sum_{k=0}^{T}\gamma_k\beta_k}\to 0, \qquad
\frac{\sum_{k=0}^{T}\gamma_k/\beta_k}{\sum_{k=0}^{T}\gamma_k\beta_k}\to 0.
\]
\end{theorem}

\begin{proof}
First, we observe the inequalities
\begin{align}
    \sum_{k=0}^T c_k &= \Psi_\la(\UU^0) - \Psi_\la(\UU^{T+1}) + \frac{1}{2 \la} \sum_{k=0}^T \ga_k^2 ( \beta_k L_\Psi+ L_\phi)^2 \nonumber\\
    \quad & \leq \Psi_\la(\UU^0) - \min_{\St} \Psi_\la + \frac{1}{2 \la} ( L_\Psi+ L_\phi/\beta_0)^2 \sum_{k=0}^T (\ga_k \beta_k)^2,\label{eq:ck-est}
\end{align}
\begin{align}
    \sum_{k=0}^T \frac{b_k^2}{\abs{a_k}} &=  \frac{2 L_\phi^2}{\la}  \sum_{k=0}^T \ga_k \left(\beta_k(\frac{1}{\la} - L_{\nabla \Psi} - L_\Psi  ) -L_\phi \right)^{-1} \leq \frac{2 L_\phi^2}{\la}  L_\la^{-1} \sum_{k=0}^T \ga_k/\beta_k, \label{eq:bk2ak-est}
\end{align}
\begin{align}
    \sum_{k=0}^T \abs{a_k} &=\frac{1}{2 \la} \sum_{k=0}^T \ga_k \left( \beta_k(\frac{1}{\la} - L_{\nabla \Psi} - L_\Psi  ) -L_\phi\right) \geq \frac{1}{2 \la} L_\la \sum_{k=0}^T \ga_k \beta_k ,\label{eq:ak-est}
\end{align}
With these and Proposition \ref{prop:Theta_esimate} we arrive at
\begin{align*}
\left(\sum_{k=0}^T \abs{a_k} \right) \underset{0\leq k\leq T}{\min} \Theta^2_{\la \Psi}(\mathbf U^k) & \leq  \sum_{k=0}^T \abs{a_k} \Theta^2_{\la \Psi}(\mathbf U^k) \leq \frac{1}{\la^2} \left( \sum_{k=0}^T \frac{b_k^2}{\abs{a_k}} + 2 \sum_{k=0}^T c_k\right) \\
&\leq \frac{2}{\la^2}(\Psi_\la(\UU^0) - \min_{\St} \Psi_\la) + \frac{1}{\la^3} ( L_\Psi+ L_\phi/\beta_0)^2 \sum_{k=0}^T (\ga_k \beta_k)^2 \\
&+ \frac{2 L_\phi^2}{\la^3}  L_\la^{-1} \sum_{k=0}^T \ga_k/\beta_k.
\end{align*}
Moreover, note that $\Psi_{\lambda}(\UU^{0})\leq \Psi(\UU^{0})$, and $\Psi_{\la}(\UU)\geq \min_{\bV}\Psi(\bV)=\min_{\St}\Psi$. Hence, 
$$
\Psi_{\la}(\UU^{0})-\min_{\St}\Psi_{\la}\leq \Psi(\UU^{0})-\min_{\St}\Psi.
$$
Putting $\sum_{k=0}^T \abs{a_k}$ on the other side and multiplying the numerator and denominator of the fraction by $\la L_\la$ yields the result.
\end{proof}

Theorem~\ref{conv_stationary} provides a general upper bound on the stationarity measure that accounts for the feasibility violation. It also specifies conditions on the step sizes and penalty parameters that guarantee convergence to a stationary point of the lower-level problem. The bound depends on the auxiliary parameter $\lambda>0$, which must be chosen sufficiently small and is then kept fixed. As $\lambda\to0$, we have $\lambda L_\lambda\to1$, while $L_\lambda$, which appears as the coefficient of the quadratic term, becomes arbitrarily large. Conversely, as $\lambda$ approaches its upper admissible bound, the constant in the denominator becomes arbitrarily small, causing the overall bound to deteriorate. 

In the following two corollaries, we specify the step sizes and penalty parameters and show that, under the differentiability of $\Psi$, the convergence rate of the stationarity measure can be transferred to the lower-level stationarity measure. In particular, we obtain an explicit convergence rate for $\|\nabla_{\St}\Psi(\UU^k)\|$,
rather than only for the Moreau-envelope-based stationarity measure $\Theta_{\la\Psi}(\UU^k)$.
The first choice of step sizes and penalty parameters yields a convergence rate for the lower-level stationarity and guarantees the existence of a subsequence converging to a stationary point. The second corollary employs slightly more restrictive step sizes, resulting in a marginally worse convergence rate, but guarantees that every limit point of a convergent subsequence is a stationary point of the lower-level problem.

\begin{corollary}\label{cor:subseq}
Let $\{\mathbf U^k\}$ be the sequence generated by \eqref{iterates} with arbitrary initialization. Let $\beta_k \rightarrow \infty$ be increase monotonically such that $\sum_{k=0}^T 1/(\beta_k^2 \sqrt{k+1}) \in o(\sqrt{T})$ and set $\ga_k = \tilde\ga/(\beta_k \sqrt{k+1})$ for some $\tilde \ga >0$. Then, for every $\lambda < (L_{\nabla \Psi} + L_\Psi + L_\phi/\beta_0)^{-1}$, we have the rate
\begin{align}
    \underset{0\leq k\leq T}{\min} \,\norm{\nabla_{\St} \Psi( \UU^k)} \in \mathcal{O}\left(\frac{\sqrt{\log{T}}}{T^{1/4}}\right)
\end{align}
and the sequence $\{\UU^k\}$ admits a convergent subsequence $\UU^{\tau(k)} \rightarrow \UU^*$ with $\nabla_{\St} \Psi(\UU^*) = 0$.
\end{corollary}

\begin{proof}
It follows from \eqref{eq:riemann} that $\snorm{\nabla_{\St} \Psi( \UU^k)} \leq (L_{\nabla_{\St}\Psi} \la +1) \snorm{\nabla_{\St} \Psi(\prox_{\la \Psi}( \UU^k))}$.
Inserting $\ga_k$ and $\beta_k$ in Theorem~\ref{conv_stationary} and noting that $\sum_{k=0}^T \frac{1}{\sqrt{k+1}} > \sqrt{T+1}$ and $\sum_{k=0}^T \frac{1}{k+1} < \log(T+1)+1$ gives the rate. 
Moreover, Theorem \ref{conv_stationary} implies that there is a subsequence $\{\UU^{k_n}\}$ such that $\snorm{\nabla_{\St} \Psi ( \prox_{\la \Psi}(\UU^{k_n}))}_F \leq \Theta(\UU^{k_n}) \rightarrow 0 $. Since the Stiefel manifold is compact, we can pass to another convergent subsequence $\UU^{k_{n_i}} \rightarrow \UU^* \in \St(K,R)$.
 Then it follows that $\prox_{\la \Psi}(\UU^{k_{n_i}}) \rightarrow \UU^*$ and the continuity of the gradient implies
 \[ \norm{\nabla_{\St} \Psi ( \prox_{\la \Psi}(\UU^{k_{n_i}}))}_F \rightarrow \norm{\nabla_{\St} \Psi (\UU^*)}_F  = 0. \hfill \text{\qedhere}\]
\end{proof}

\begin{corollary}\label{prop_conv}
Let $\{\mathbf U^k\}$ be the sequence generated by \eqref{iterates} with arbitrary initialization. Let $\beta_k \rightarrow \infty$ be monotonically increasing with $\beta_k\in \mathcal{O}(k^{1/4})$ and set $\ga_k= \tilde\ga/(\beta_k (k+1)^{1/2+\eps})$ for some fixed $\eps, \tilde \ga > 0$. Then, for every $\lambda < (L_{\nabla \Psi} + L_\Psi + L_\phi/\beta_0)^{-1}$, we have the rate
\begin{equation}
    \underset{0\leq k\leq T}{\min} \norm{\nabla_{\St} \Psi( \UU^k)}^2  \leq \frac{1}{T+1}\sum_{k=0}^T \norm{\nabla_{\St} \Psi( \prox_{\la \Psi}(\UU^k))}^2 \in \mathcal{O} (T^{-1/2 +\eps}).
\end{equation}
Especially, the limit point $\UU^*$ of any convergent subsequence $\UU^{\tau(k)}$ satisfies $\nabla_{\St} \Psi(\UU^*) = 0$.
\end{corollary}

\begin{proof}
It follows from \eqref{eq:riemann} that $\snorm{\nabla_{\St} \Psi( \UU^k)} \leq (L_{\nabla_{\St}\Psi} \la +1) \snorm{\nabla_{\St} \Psi(\prox_{\la \Psi}( \UU^k))}$.
By Proposition~\ref{prop:Theta_esimate} we have
    \begin{align*}
        \frac{1}{T+1}\sum_{k=0}^T \Theta^2_{\la \Psi}(\mathbf U^k) 
        &\leq \frac{1}{\la^2 (T+1)}\sum_{k=0}^T \frac{1}{\abs{a_k}} \left(\frac{b_k^2}{\abs{a_k}} + 2c_k\right) \leq \frac{1}{\la^2 (T+1)} \, \underset{0\leq k\leq T}{\max} \,\left(\frac{1}{\abs{a_k}} \right)\cdot \sum_{k=0}^T  \left(\frac{b_k^2}{\abs{a_k}} + 2c_k\right) 
    \end{align*}
Further, $\underset{0\leq k\leq T}{\max} \,\frac{1}{\abs{a_k}} \in \mathcal{O}(1/(\ga_T \beta_T)) = \mathcal{O}(T^{1/2+\eps})$, and by~\eqref{eq:bk2ak-est} and~\eqref{eq:ck-est}
\begin{align*}
    \sum_{k=0}^T  \frac{b_k^2}{\abs{a_k}} + 2c_k 
    &\leq\frac{2 L_\phi^2}{\la}  L_\la^{-1} \sum_{k=0}^T \ga_k/\beta_k  + 2( \Psi_\la(\UU^0) - \min \Psi_\la ) + \frac{1}{2 \la} ( L_\Psi+ L_\phi/\beta_0)^2 \sum_{k=0}^T (\ga_k \beta_k)^2\\
    &\leq \frac{2 L_\phi^2}{\la}  L_\la^{-1} \sum_{k=0}^T \frac{\tilde \ga}{(k+1)^{1+\eps}}  + 2\Delta_{\Psi}^{0}+ \frac{1}{2 \la} ( L_\Psi+ L_\phi/\beta_0)^2 \sum_{k=0}^T \frac{\tilde \ga^2}{(k+1)^{1+2\eps}}\\
    & \in \mathcal{O}(1)
\end{align*}
Since the average converges to zero, the entire sequence satisfies $\snorm{\nabla_{\St}\Psi\left(\prox_{\lambda\Psi}(\UU^k)\right)}_F \to 0$.
The remainder of the proof follows in the same manner as in Corollary~\ref{cor:subseq}.
\end{proof}

The corollaries specify the step sizes up to a multiplicative constant, while the penalty parameters still allow some flexibility. We briefly provide examples of penalty parameter sequences that satisfy the assumptions of the analysis. 
\begin{enumerate}
\item \textbf{Polynomial growth:} Let
\[
\beta_k = C (k+1)^p,
\qquad C>0, p>0.
\]
Then
\[
\sum_{k=0}^T \frac{1}{\beta_k^2 \sqrt{k+1}} = \sum_{k=0}^T \frac{1}{(k+1)^{1/2 + 2 p}} \leq \begin{cases}
\frac{T^{1/2-2 p}-1}{1/2-2 p} \in o(\sqrt{T}) \quad &p<\frac{1}{4}\\
\log(T) \in o(\sqrt{T}) \quad &p=\frac{1}{4}\\
\in \mathcal{O}(1) \subset o(\sqrt{T}) \quad & p > \frac{1}{4}.
\end{cases}
\]
Hence, the assumptions of the first corollary are satisfied for all $p>0$, whereas the second corollary requires $p \leq 1/4$. 
\item \textbf{Logarithmic growth:} Let
\[
\beta_k \ge C \sqrt{\log(k+2)}, \qquad C >0.
\]
Then
\[
\frac{1}{\sqrt{T}}
\sum_{k=0}^T
\frac{1}{\sqrt{k+1}\beta_k^2}
\in 
\mathcal{O}\!\left(\frac{\log\log T}{\sqrt{T}}\right)
\to 0.
\]
To satisfy the assumptions of the second corollary, we may choose $\beta_k=C\sqrt{\log(k+2)}$ since $\sqrt{\log(k+2)}\in\mathcal{O}(k^{1/4})$.
\end{enumerate}

Next, we analyze upper-level stationarity. While for the lower level we track the evolution of the auxiliary sequence $\{\prox_{\lambda\Psi}(\UU^k)\}_{k}$, for the upper level we instead consider the sequence $\hat \UU^k \coloneqq \prox_{\la \Omega_k/\beta_k}(\UU^k)$ and derive a general upper bound for the time-dependent stationarity measure $\beta_k\Theta_{\lambda\Omega_k/\beta_k}$. This measure provides a bound on the distance to the subdifferential $\partial_{\St} \Omega_k(\hat{\UU}^k)$ to zero. The resulting bound directly yields conditions on the penalty parameters under which a subsequence of these subgradients converges to zero. 
In contrast to the lower-level case, however, this does not imply that the subgradients of $\Omega_k$ evaluated along the algorithm-generated sequence $\{\UU^k\}$ converge to zero. This is due to the lack of differentiability of the upper-level objective and the scaling induced by the increasing penalty parameter. Therefore, we do not aim to establish convergence of these subgradients to zero.
Instead, combined with additional regularity assumptions on $\Psi$ and suitable regularity conditions at the limit points of the algorithm, we can show that a limit point $\UU^*$ not only belongs to the lower-level stationary set
\[
N_\Psi \coloneqq \{\UU\mid \nabla_{\St}\Psi(\UU)=0\},
\]
but is, in fact, a stationary point of the upper-level problem $\phi$ when restricted to $N_\Psi$.

Before proving the bounds for the stationarity measure, we establish several technical inequalities. The following lemma shows that the auxiliary parameter $\lambda$ must be chosen sufficiently small to exploit the additional convexity of the Moreau envelope. Moreover, the proof builds on arguments from simple bilevel optimization in Euclidean spaces; however, some of these arguments fail here due to the lack of convexity of $\phi$ when restricted to the Stiefel manifold. We overcome this issue by introducing an additional auxiliary sequence along the Euclidean line segments connecting the two sequences on the manifold.

\begin{lemma}\label{lem:ineq}
Let $\{\mathbf U^k\}$ be the sequence generated by \eqref{iterates} with arbitrary initialization, step sizes $(\ga_k)_k$ and monotonically increasing penalty parameters $(\beta_k)_k$. Define $\rho \coloneqq L_{\nabla \Psi} + 4 L_\Psi + 4 L_\phi/\beta_0$ and choose $\la < 1/\rho$. Moreover, define $\hat\UU^k \coloneqq \prox_{(\la/\beta_k) \Omega_k}(\UU^k)$. For arbitrary $g^k \in \partial_{\St} \phi(\UU^k)$ and  $\bV \in \St$, the following inequalities hold:
    \begin{align}
        \left(\frac{1}{\la } - \rho\right) \norm{\UU^k - \hat \UU^k}_{F}^2 
        &\leq 8 \inner{\frac{1}{\beta_k}g^k + \nabla_{\St} \Psi(\UU^k)}{\UU^k - \hat \UU^k} \label{ineq_conv} \\
        \ga_k \beta_k \inner{\frac{1}{\beta_k}g^k + \nabla_{\St} \Psi(\UU^k)}{\UU^k - \bV} 
        &\leq \frac{1}{2} \norm{\UU^k - \bV}^2_{F} - \frac{1}{2} \norm{\UU^{k+1} - \bV}^2_{F} + \frac{1}{2}(\ga_k \beta_k)^2 \left(\frac{1}{\beta_0}L_\phi + L_\Psi\right)^2 \label{ineq_iterate}\\
        \frac{1}{2} \norm{\UU^k - \hat \UU^k}^2_{F} - \frac{1}{2} \norm{\UU^{k+1} - \hat \UU^k}^2_{F}
        &\leq \la\left[ (\beta_k^{-1}\Omega_k)_\la (\UU^k) - (\beta_{k+1}^{-1} \Omega_{k+1})_\la (\UU^{k+1}) - (\beta_{k}^{-1}-\beta_{k+1}^{-1})\min_{\UU \in \St} \phi(\UU)\right] \label{ineq_enve}
        \end{align}
\end{lemma}
\begin{proof}
For \eqref{ineq_conv} note that $\frac{1}{\beta_k} \Omega_k$ is $L_{\nabla \Psi}$-weakly convex on $\mathcal{U}$, therefore $f_k(\UU) \coloneqq \frac{1}{\beta_k} \Omega_k(\UU) + \frac{1}{2 \la } \norm{\UU - \UU^k}^2$ is $\left(1/\la - L_{\nabla\Psi}\right)$-strongly convex on $\mathcal{U}$. Define $\dot{\UU}^k \coloneqq \argmin_{\UU \in \left[\UU^k, \hat\UU^k\right]} f_k(\UU)$ where $\left[\UU^k, \hat\UU^k\right] \subset \mathcal{U}$ denotes the line segment between $\UU^k$ and $\hat \UU^k$ in the ambient Euclidean space. By the strong convexity of $f_k$ and the fact that $\left[\UU^k, \hat\UU^k\right]$ is convex, $\dot{\UU}^k$ is unique and for all $\UU \in \left[\UU^k, \hat\UU^k\right]$ holds
\[ f_k(\UU) - f_k(\dot{\UU}^k) \geq \frac{1}{2} \left(\frac{1}{\la} - L_{\nabla\Psi}\right) \norm{\UU - \dot{\UU}^k}^2_{F}.\]
Then, we have
\begin{align*}
\frac{1}{2} \left(\frac{1}{\la} - L_{\nabla\Psi}\right) \norm{\UU^k - \hat \UU^k}^2_{F} 
 &\leq \frac{1}{2} \left(\frac{1}{\la} - L_{\nabla\Psi}\right) \left( 2 \norm{\UU^k - \dot{\UU}^k}^2_{F}  + 2  \norm{\hat\UU^k - \dot{\UU}^k}^2_{F} \right) \\
 & \leq 2 \left(f_k(\UU^k) - f_k(\dot{\UU}^k) + f_k(\hat \UU^k) - f_k(\dot{\UU}^k) \right) \\
 & \leq 4 \left(f_k(\UU^k) - f_k(\dot{\UU}^k) \right),
\end{align*}
where, in the last inequality, we have used the fact that $\hat \UU^k \in \argmin_{\UU \in \St} f_k(\UU)$ and therefore $f_k(\hat \UU^k) \leq f_k(\UU^k)$. Defining $\tilde \rho = L_{\nabla\Psi} + L_\Psi + L_\phi/\beta_0 < \rho$ and using the definition of $f_k$ yields
\begin{align*}
 f_k(\UU) - f_k(\dot{\UU}^k) 
 & = \frac{1}{\beta_k}\Omega_k(\UU^k) + \frac{1}{2 \la} \norm{\UU^k - \UU^k}^2_{F} - \frac{1}{\beta_k}\Omega_k(\dot{\UU}^k) - \frac{1}{2 \la} \norm{\dot{\UU}^k - \UU^k}^2_{F} \\
 & + \frac{1}{2} \left(\frac{1}{\la} - \tilde\rho\right) \norm{\UU^k - \dot{\UU}^k}^2_{F} - \frac{1}{2} \left(\frac{1}{\la} - \tilde\rho\right) \norm{\UU^k - \dot{\UU}^k}^2_{F} \\
 & = \frac{1}{\beta_k}\Omega_k(\UU^k) - \frac{1}{\beta_k}\Omega_k(\dot{\UU}^k) - \frac{\tilde\rho}{2} \norm{\dot{\UU}^k - \UU^k}^2_{F}  - \frac{1}{2} \left(\frac{1}{\la} - \tilde\rho\right) \norm{\UU^k - \dot{\UU}^k}^2_{F} \\
 & \leq \frac{1}{\beta_k}\Omega_k(\UU^k) - \frac{1}{\beta_k}\Omega_k(\dot{\UU}^k) - \frac{L_{\nabla\Psi} + L_\Psi + L_\phi/\beta_0}{2} \norm{\dot{\UU}^k - \UU^k}^2_{F} \\
 & \leq \inner{\frac{1}{\beta_k}\omega^k}{\UU^k - \dot{\UU}^k} - \frac{L_\Psi + L_\phi/\beta_0}{2} \norm{\dot{\UU}^k - \UU^k}^2_{F} , 
\end{align*}
where the last estimate is a consequence of $\frac{1}{\beta_k}\Omega_k$ being $L_{\nabla\Psi}$-weakly convex on $\mathcal{U}$ and   $\omega^k \in \partial \Omega_k(\UU^k)$ being an Euclidean subgradient.
Now observe that $\dot{\UU}^k = t \UU^k + (1-t) \hat \UU^k$ for some $t \in [0,1]$ and hence 
\[ \UU^k - \dot{\UU}^k = (1-t) \UU^k - (1-t)\hat \UU^k = (1-t) (\UU^k - \hat \UU^k). \]
From Theorem \ref{thm:subgrad_ineq} we deduce the bound 
\[\inner{\frac{1}{\beta_k}\omega^k}{\UU^k - \hat \UU^k} \leq \inner{\frac{1}{\beta_k} \proj_{T_{\UU^k}\St}(\omega^k)}{\UU^k - \hat \UU^k} + \frac{L_\Psi + L_\phi/\beta_0}{2}\norm{\hat \UU^k - \UU^k}^2_{F}.
\]
Combined with non-negativity of the right hand side, this leads to
 \begin{align*}
 \inner{\frac{1}{\beta_k}\omega^k}{\UU^k - \dot{\UU}^k} 
 & - \frac{L_\Psi + L_\phi/\beta_0}{2} \norm{\dot{\UU}^k - \UU^k}^2 \\
 & =  (1-t)\left( \inner{\frac{1}{\beta_k}\omega^k}{\UU^k - \hat \UU^k} - (1-t) \frac{L_\Psi + L_\phi/\beta_0}{2} \norm{\UU^k - \hat \UU^k}^2_{F}\right)\\
 & \leq (1-t) \left(\inner{\frac{1}{\beta_k}g^k + \nabla_{\St}\Psi(\UU^k)}{\UU^k - \hat \UU^k} + t \frac{L_\Psi + L_\phi/\beta_0}{2} \norm{\UU^k - \hat \UU^k}^2_{F} \right) \\
 & \leq \inner{\frac{1}{\beta_k}g^k + \nabla_{\St}\Psi(\UU^k)}{\UU^k - \hat \UU^k} + \frac{L_\Psi + L_\phi/\beta_0}{2} \norm{\UU^k - \hat \UU^k}^2_{F}.
 \end{align*}
Putting all together, we conclude with
\begin{align*}
    \frac{1}{8} \left(\frac{1}{\la} - L_{\nabla\Psi}\right) \norm{\UU^k - \hat \UU^k}^2_{F} 
    & \leq f_k(\UU^k) - f_k(\dot{\UU}^k) \\
    & \leq \inner{\frac{1}{\beta_k}g^k + \nabla_{\St} \Psi(\UU^k)}{\UU^k - \hat \UU^k} + \frac{L_\Psi + L_\phi/\beta_0}{2} \norm{\UU^k - \hat \UU^k}^2_{F}.
\end{align*}
Finally, subtracting $\frac{L_\Psi + L_\phi/\beta_0}{2} \norm{\UU^k - \hat \UU^k}^2_{F}$ from both sides and using the definition of $\rho$ yields \eqref{ineq_conv}.

Next, we consider
\begin{align*}
    \norm{\UU^{k+1} - \bV}^2_{F} 
    & = \norm{ \proj_{\St}\Bigl(\UU^k -\ga_k \bigl( g^k+\beta_k\nabla_{\St} \Psi(\UU^k) \bigr) \Bigr) - \bV}^2_{F} \\
    &\leq \norm{ \UU^k -\ga_k \bigl( g^k+\beta_k\nabla_{\St} \Psi(\UU^k) \bigr) - \bV}^2_{F} \\
    & = \norm{\UU^k - \bV}^2_{F} - 2 \ga_k \inner{g^k + \beta_k \nabla_{\St} \Psi(\UU^k)}{\UU^k- \bV} + \ga_k^2 \norm{g^k + \beta_k \nabla_{\St} \Psi(\UU^k)}^2_{F} 
\end{align*}
where $\bV \in \St(K,R)$ and Lemma \ref{polar_retraction} have been used. To bound the last quadratic term on the right-hand side, we observe
\[ \norm{g^k + \beta_k \nabla_{\St} \Psi(\UU^k)}^2_{F} = \beta_k^2 \norm{\frac{1}{\beta_k}g^k + \nabla_{\St} \Psi(\UU^k)}^2_{F} \leq \beta_k^2 \left( \frac{1}{\beta_0} L_\phi + L_\Psi \right)^2.
\]
Substituting this into the previous display and rearranging the terms yields \eqref{ineq_iterate}.

For the final relation, we start with 
\[ (\beta_{k+1}^{-1} \Omega_{k+1})_\la (\UU^{k+1}) = \min_{\UU \in \St} \left\{\beta_{k+1}^{-1} \Omega_{k+1} (\UU) + \frac{1}{2 \la} \norm{\UU - \UU^{k+1}}^2_{F}\right\} \leq \beta_{k+1}^{-1} \Omega_{k+1}(\hat \UU^k) + \frac{1}{2 \la} \norm{\hat \UU^k - \UU^{k+1}}^2_{F}
\]
and by the definition of $\hat \UU^k$ we have 
\[ (\beta_{k}^{-1} \Omega_{k})_\la (\UU^k) =  \beta_{k}^{-1} \Omega_{k}(\hat \UU^k) + \frac{1}{2 \la} \norm{\hat \UU^k - \UU^{k}}^2_{F}.
\]
Adding these two together, we conclude with
\begin{align*}
    \frac{1}{2 \la}\norm{\hat \UU^k - \UU^k}^2_{F} 
    & - \frac{1}{2 \la } \norm{\hat \UU^k - \UU^{k+1}}^2_{F} \\
    &\leq \left(  (\beta_k^{-1}\Omega_k)_\la (\UU^k)- \beta_k^{-1} \Omega_k(\hat \UU^k) \right) + \left(  \beta_{k+1}^{-1}\Omega_{k+1} (\hat \UU^k)- (\beta_{k+1}^{-1} \Omega_{k+1})_\la (\UU^{k+1}) \right) \\
    & \leq (\beta_k^{-1}\Omega_k)_\la (\UU^k) - (\beta_{k+1}^{-1} \Omega_{k+1})_\la (\UU^{k+1}) - (\beta_{k}^{-1}-\beta_{k+1}^{-1})\phi(\hat \UU^k).
\end{align*}
By monotonicity $\beta_{k}^{-1}-\beta_{k+1}^{-1} \geq 0$ and estimating $\phi(\hat \UU^k)$ with its minimum, we arrive at \eqref{ineq_enve}. 
\end{proof}

\begin{theorem} \label{thm:prox_conv}
Let $\{\mathbf U^k\}$ be the sequence generated by \eqref{iterates} with arbitrary initialization, step sizes $(\ga_k)_k$ and monotonically increasing penalty parameters $(\beta_k)_k$. Define $\rho \coloneqq L_{\nabla \Psi} + 4 L_\Psi + 4 L_\phi/\beta_0$.  Then, for every $\lambda < 1/\rho$, the following estimate holds
\begin{equation}\label{eq:prox_conv}
    \min_{0 \leq k \leq T} \beta_k^2  \Theta_{\la \Omega_k/\beta_k}(\mathbf{U}^k) ^2\leq \frac{8}{1- \la \rho} \frac{\beta_0 \Omega_0(\UU^0) + \beta_{T+1}^2 C_{\max} + 2 \beta_{T+1} \abs{\phi^*} + \frac{L^2}{2 \la}\sum_{k=0}^T \beta_k^2 (\ga_k \beta_k)^2 }{ \sum_{k=0}^{T} \ga_k \beta_k},
\end{equation}
where $C_{\max} \coloneqq \beta_0^{-1} \max_{\St} \phi + \max_{\St} \Psi + \sabs{\min_{\St} \Psi}$, $\phi^* \coloneqq \min_{\St} \phi$ and $L \coloneqq \beta_0^{-1}L_\phi + L_\Psi$. In particular, the right-hand side of \eqref{eq:prox_conv} converges to zero as $T \to \infty$ if and only if the following design criteria for the step-size and penalty sequences are satisfied:
\[
\sum_{k=0}^{T}\gamma_k\beta_k \to\infty, \qquad
\frac{\sum_{k=0}^{T}\beta_k^2(\gamma_k\beta_k)^2}{\sum_{k=0}^{T}\gamma_k\beta_k}\to 0, \qquad
\frac{\beta_{T+1}^2}{\sum_{k=0}^{T}\gamma_k\beta_k}\to 0.
\]
\end{theorem}

\begin{proof}
Define $\hat\UU^k \coloneqq \prox_{(\la/\beta_k) \Omega_k}(\UU^k)$ and as defined in \eqref{eq:theta} we have 
\[
\Theta_{(\la/\beta_k)\Omega_k}(\mathbf{U}^k) = \lambda^{-1}\| \hat\UU^k-\UU^k\|.
\]
The inequalities of Lemma~\ref{lem:ineq} yield the following estimate:
\begin{align*}
     (\ga_k \beta_k)  \beta_k^2  \la^{-2} \norm{\UU^k- \hat \UU^k}^2 
     & \leq \frac{8}{\la(1- \la \rho)}  \beta_k^2 (\ga_k \beta_k) \inner{\frac{1}{\beta_k}g^k + \nabla_{\St} \Psi(\UU^k)}{\UU^k - \hat \UU^k}\\
     & \leq \frac{8}{(1- \la \rho)} \beta_k^2 \big[ (\beta_k^{-1}\Omega_k)_\la (\UU^k) - (\beta_{k+1}^{-1} \Omega_{k+1})_\la (\UU^{k+1}) \\
     & - (\beta_{k}^{-1}-\beta_{k+1}^{-1}) \phi^* +  \frac{1}{2 \la}(\ga_k \beta_k)^2 \left(\frac{1}{\beta_0}L_\phi + L_\Psi\right)^2 \big].
\end{align*}

To obtain a recursion, consider that $(\beta_k)_k$ is monotonically increasing, we estimate
\[ - \beta_k^2 (\beta_{k}^{-1}-\beta_{k+1}^{-1})\phi^* \leq 
\beta_k \beta_{k+1} (\beta_{k}^{-1}-\beta_{k+1}^{-1})\abs{\phi^*} \leq (\beta_{k+1} - \beta_k) \abs{\phi^*}.
\]
Additionally, we have 
\begin{align*}
    - \beta_k^2 (\beta_{k+1}^{-1} \Omega_{k+1})_\la(\UU^{k+1}) 
    &  = - \beta_{k+1}^2 (\beta_{k+1}^{-1} \Omega_{k+1})_\la(\UU^{k+1}) + (\beta_{k+1}^2- \beta_k^2) (\beta_{k+1}^{-1} \Omega_{k+1})_\la(\UU^{k+1}) \\
    & \leq - \beta_{k+1}^2 (\beta_{k+1}^{-1} \Omega_{k+1})_\la(\UU^{k+1}) + (\beta_{k+1}^2- \beta_k^2) \left(\frac{1}{\beta_0} \max_{\UU \in \St} \phi(\UU) + \max_{\UU \in \St} \Psi(\UU) \right)
\end{align*}

We use the upper bound 
\[\left(\frac{1}{\beta_{k+1}}\Omega_{k+1}\right)_{\la}(\UU^{k+1}) \leq \frac{1}{\beta_{k+1}}\Omega_{k+1}(\UU^{k+1}) \leq \frac{1}{\beta_0} \max_{\UU \in \St} \phi(\UU) + \max_{\UU \in \St} \Psi(\UU) \]
in order to obtain the expression $- \beta_{k+1}^2 (\beta_{k+1}^{-1} \Omega_{k+1})_\la(\UU^{k+1})$ for the recursion.
Summing over $k=0, \ldots, T$ yields
\begin{align*}
    \sum_{k=0}^{T} (\ga_k \beta_k) \beta_k^2 \Theta_{(\la/\beta_k) \Omega_k} (\mathbf{U}^k)^2
    & \leq \frac{8}{1- \la \rho} \Big[  \beta_0^2 (\beta_0^{-1}\Omega_0)_\la (\UU^0) - \beta_{T+1}^2(\beta_{T+1}^{-1} \Omega_{T+1})_\la (\UU^{T+1})\\
    & + (\beta_{T+1}^2-\beta_0^2) \frac{1}{\beta_0} \max_{\St} \phi + \max_{\St} \Psi\\
     & + (\beta_{T+1}-\beta_{0}) \abs{\phi^*} +  \frac{1}{2 \la} \left(\frac{1}{\beta_0}L_\phi + L_\Psi\right)^2 \sum_{k=0}^{T} \beta_k^2 (\ga_k \beta_k)^2 \Big]
\end{align*}
Moreover, note that $\beta_0^2(\beta_0 \Omega_0)_{\la}(\UU^0) \leq \beta_0 \Omega_0(\UU^0)$  and $\beta_{T+1}^{-1} \Omega_{T+1})_\la (\UU^{T+1}) \geq - ( \beta_{T+1}^{-1} \sabs{\phi^\ast} + \sabs{\min_{\St}\Psi})$. Hence, 
$$
\beta_0^2 (\beta_0^{-1}\Omega_0)_\la (\UU^0) - \beta_{T+1}^2(\beta_{T+1}^{-1} \Omega_{T+1})_\la (\UU^{T+1})\leq \beta_0 \Omega_0(\UU^0) + \beta_{T+1} \sabs{\phi^\ast} + \beta_{T+1}^2 \sabs{\min_{\St}\Psi}.
$$
By positivity, we have $\min_{0\leq k\leq T}\beta_k^2 \Theta_{(\la/\beta_k) \Omega_k}(\mathbf{U}^k)^2\sum_{k=0}^{T}\gamma_{k}\beta_{k}\leq \sum_{k=0}^{T} \ga_k \beta_k \beta_k^2 \Theta_{\la/\beta_k \Omega_k} (\mathbf{U}^k)^2$. Canceling the negative terms yield the bound \eqref{eq:prox_conv}.
\end{proof}

Theorem~\ref{thm:prox_conv} provides a general upper bound on the stationarity measure $\beta_k\Theta_{\lambda\Omega_k/\beta_k}$ that accounts for the penalized objective function. It also specifies conditions on the step sizes and penalty parameters that guarantee convergence to a stationary point. The dependence on the auxiliary parameter $\lambda>0$ is similar to that in Theorem~\ref{conv_stationary}: as $\lambda\to0$, the coefficient of the quadratic term becomes arbitrarily large, while as $\lambda$ approaches its upper admissible bound, the constant in the denominator becomes arbitrarily small, causing the overall bound to deteriorate.

In the following corollary, we specify the step sizes and penalty parameters to obtain a convergence rate for the penalized problem.

\begin{corollary}\label{cor:comb_conv}
Let $\{\mathbf U^k\}$ be the sequence generated by \eqref{iterates} with arbitrary initialization. Set $\ga_k = \tilde \ga/(\beta_k \sqrt{k+1})$ with $\beta_k = C(k+1)^p$ for some $C>0$ and $p \in(0,1/2)$. Define $\hat\UU^k = \prox_{\la\Omega_k/\beta_k}(\UU^k)$. Then
\begin{equation}
\min_{0 \leq k \leq T} \dist(0,\partial_{\St}\Omega_{k}(\hat{\UU}^{k}))
\in \mathcal{O}\left(T^{-(1/4-p)}\right).
\end{equation}
In particular, the stationarity measure converges to zero as $T\to\infty$ for any $p<1/4$.
\end{corollary}
\begin{proof}
Recall from~(\ref{eq:theta}) that $\hat{\UU}^{k}$ fulfills
$\dist(0,\partial_{\St}\Omega_{k}(\hat{\UU}^{k}))\leq\beta_{k}\Theta_{\lambda\Omega_{k}/\beta_k}(\UU^k)$. Using the definition of the penalty and step-size parameter, we obtain for $p <1/2$
\begin{align*}
\sum_{k=0}^{T}\gamma_{k}\beta_{k}&=\tilde{\gamma}C\sum_{k=0}^{T}\frac{1}{\sqrt{k+1}}\geq\tilde{\gamma}C\sqrt{T+1} \\
\sum_{k=0}^{T}\gamma_{k}^{2}\beta^{4}_{k}&=\tilde{\gamma}^{2}C^4\sum_{k=0}^{T}\frac{1}{(k+1)^{1-2p}}=\tilde{\gamma}^{2}C^4 \left(1+\sum_{k=1}^{T}\frac{1}{(k+1)^{1-2p}}\right)\\
&\leq \tilde{\gamma}^{2}C^4 \left(1+\int_{0}^{T}\frac{1}{(x+1)^{1-2p}}\dif x \right)\leq \tilde{\gamma}^{2}C^4 \left(1 +\frac{(T+1)^{2p}}{2p}\right).
\end{align*}
From Theorem \ref{thm:prox_conv}, we deduce that 
\begin{align*}
    \min_{0\leq k\leq T}&\dist(0,\partial_{\St}\Omega_{k}(\hat{\UU}^{k})) \\
    &\leq \sqrt{\frac{8}{1-\lambda\rho}\frac{\beta_0 \Omega_0(\UU^0) + C^2(T+2)^{2p} C_{\max} + 2 C(T+2)^{p} \abs{\phi^{*}}+\frac{L^{2}\tilde{\gamma}^{2}C^4}{2\lambda}\left(1+\frac{(T+1)^{2p}}{2p}\right)}{\tilde{\gamma}C\sqrt{T+1}}}\\
    &\in \scrO(T^{p-1/4}). && \text{\qedhere}
\end{align*}
\end{proof}

The corollary shows that the distance to the Riemannian subdifferential can be made arbitrarily small. However, this does not yet imply that the corresponding Riemannian subgradients converge. In contrast to the lower-level problem, the subdifferential here is set-valued rather than given by a continuous gradient. Moreover, since $\phi$ is assumed to be convex only in the ambient Euclidean space, the corresponding outer semicontinuity of the Riemannian subdifferential is not immediate. The following lemma shows that, nevertheless, convergent subsequences of Riemannian subgradients have limits that remain in the Riemannian subdifferential.

\begin{lemma}\label{conv_subgradient}
Let $\{\bV^k\}$ be a convergent sequence in $\St(K,R)$ with limit point $\bV^*$ and $g^k \in \partial_{\St} \phi(\bV^k)$ be a collection of Riemannian subgradients of the convex function $\phi$. Then  $\{g^k\}$ admits a converging subsequence $g^{\tau(k)} \rightarrow g^* \in \partial_{\St} \phi(\bV^*)$.
\end{lemma}

\begin{proof}
Since $\St(K,R)$ is compact and $\phi$ is convex, the Euclidean subgradients of $\phi$ are uniformly bounded on $\St(K,R)$, i.e.
\[(\exists M > 0)(\forall \bV \in \St(K,R)) \colon \forall \tilde g \in \partial \phi(\bV) \colon \quad \norm{\tilde g} \leq M.\]
Hence, every collection of Euclidean subgradients $\{\tilde g^k \in \partial \phi(\bV^k)\}$ is bounded.
Since we are in finite dimensional spaces, every bounded sequence admits
a convergent subsequence. Thus, we can pick a convergence subsequence 
$\tilde g^{\tau(k)} \to \tilde g^*$.
We show that $\tilde g^* \in \partial_{\St} \phi(\bV^*)$. For every $k$ and every $\bV \in \mathbb R^{K\times R}$, 
\[\phi(\bV)\ge\phi(\bV^{\tau(k)})+\langle \tilde{g}^{\tau(k)}, \bV-\bV^{\tau(k)} \rangle .\]
Passing to the limit as $k\to\infty$ and using the continuity of convex
functions yields $\tilde g^* \in \partial \phi(\bV^*)$. Now we show the same for the Riemannian subgradients. We know that every Riemannian subgradient is the orthogonal projection of an Euclidean subgradient onto the tangent space (cf. Appendix \ref{app:basic}):
\[ g^k = \proj_{T_{\bV^k}\St} (\tilde g^k).\]
Since orthogonal projection is nonexpansive we have $\norm{g^k} \leq \norm{\tilde g^k}$ and hence every collection $\{g^k \in \partial_{\St} \phi(\bV^k)\}$ is bounded and admits a converging subsequence $g^{\tau(k)} \rightarrow g^*$. Then, by the continuity of the tangent-space projection 
\[g^* = \proj_{T_{\bV^*}\St}(\tilde g^*) \in \partial_{\St}\phi(\bV^*).\hfill \text{\qedhere}\]
\end{proof}
In summary, the above results show that, with $\hat\UU^k = \prox_{\lambda\Omega_k/\beta_k}(\UU^k)$, there exists a subsequence $\{\tau(k)\}$ and a limit point $\UU^*$ such that
\[
\UU^{\tau(k)}\to\UU^*,
\qquad
\hat\UU^{\tau(k)}\to\UU^*,
\qquad
\UU^* \in N_{\Psi} = \{\UU \mid \nabla_{\St}\Psi(\UU)=0\}.
\]
Moreover, there exist $\hat g^k\in\partial_{\St}\phi(\hat\UU^k)$ and $g^*\in\partial_{\St}\phi(\UU^*)$ such that
\[
\hat g^{\tau(k)}\to g^*,
\qquad
\left\|
\hat g^{\tau(k)}
+\beta_{\tau(k)}
\nabla_{\St}\Psi(\hat\UU^{\tau(k)})
\right\|\to0.
\]
These stationarity properties provide useful information about the asymptotic behavior of the algorithm. To obtain a more precise characterization of the limit points and, in particular, to relate them to the upper-level problem, we impose additional regularity assumptions on $\Psi$ and at the limit points $\UU^*$. Under these assumptions, we can characterize the limit points of the algorithm more precisely: they are, in fact, stationary points of the upper-level objective $\phi$ restricted to the stationary set of $\Psi$.


\paragraph{Taylor expansion on the Stiefel manifold.}
In the following, we use a second-order Taylor expansion of $\Psi$ on the Stiefel manifold. We assume that $\Psi$ is twice continuously differentiable. Unlike in Euclidean space, where Taylor expansions are formulated along straight lines, on a Riemannian manifold they are formulated along geodesics, which are the manifold analogues of straight lines. Moreover, tangent vectors at different points belong to different tangent spaces and therefore cannot be directly compared. Parallel transport provides a way to move tangent vectors between tangent spaces while preserving the Riemannian metric.

A detailed introduction to geodesics, the exponential and logarithmic maps, and parallel transport is provided in Appendix~\ref{app:taylor}; see also \cite[Chapter 10]{Boumal_2023}. Here, let $c(t) = \mathrm{Exp}_{\bA}(t\xi)$ be the geodesic connecting two sufficiently close points $\bA,\bB\in\St$, with $\bA=c(0)$, $\bB=c(1)=\mathrm{Exp}_{\bA}(\xi)$ and 
$\xi = \mathrm{Log}_{\bA}(\bB) \in T_{\bA}\St$. 
Thus, $\xi$ is the tangent vector at $\bA$ pointing along the geodesic from $\bA$ to $\bB$. We denote parallel transport along $c$ from $c(0)$ to $c(t)$ by
\[
P_{t\xi}=PT^c_{t\leftarrow 0}
\colon T_{c(0)}\St \rightarrow T_{c(t)}\St.
\]

Let $\{\bV^k\}$ be a convergent sequence in $\St(K,R)$ with limit point $\bV^*$ such that $\nabla_{\St} \Psi(\bV^*) = 0$. For $k \in \N$ sufficiently large, let $v_k =  \mathrm{Log}_{\bV^*}(\bV^k) \in T_{\bV^*}\St$ denote the initial velocity of the geodesic $c_k(t) = \mathrm{Exp}_{\bV^*}(t v_k)$  with $c_k(0) = \bV^*$ and $c_k(1) = \bV^{k}$. Since $\nabla_{\St}\Psi(\bV^*) = 0$, Taylor's theorem yields
\begin{equation}\label{eq:taylor}
P_{v_{k}}^{-1} \nabla_{\St} \Psi(\bV^{k}) = \int_{0}^{1} (P_{t v_{k}}^{-1}\circ \mathrm{Hess}_{\St} \Psi(c_k(t))\circ P_{t v_{k}})[v_{k}] \dif t.
\end{equation}
We define the family of linear operators $\{B^{k}\}_{k}\colon T_{\bV^*}\St \rightarrow T_{\bV^*}\St$ by
\begin{equation}\label{eq:B_k}
B^k \coloneqq \int_{0}^{1} P_{t v_{k}}^{-1}\circ \mathrm{Hess}_{\St} \Psi(c_k(t)) \circ P_{t v_{k}} \dif t.
\end{equation}
Since $c_k(t)\to\bV^*$ uniformly in $t$, 
$\mathrm{Hess}_{\St}\Psi(c_k(t))\to
\mathrm{Hess}_{\St}\Psi(\bV^*)$ uniformly in $t$. Moreover,
$P_{t v_k}\to\Id$ and $P_{t v_k}^{-1}\to\Id$ uniformly in $t$.
Consequently,
\[B^k \to \mathrm{Hess}_{\St}\Psi(\bV^*) \quad \text{as} \quad k \rightarrow \infty.\]

Throughout the remainder of this section, we impose the following assumptions.
\begin{assumption}\label{ass:Psi2}
$\Psi$ is twice continuously differentiable on $\St(K,R)$.
\end{assumption}
\begin{assumption}[Constant-rank condition]\label{ass:rank}
Let $\{\bV^k\}$ be a convergent sequence in $\St(K,R)$ with
$\bV^k\to\bV^*$. The corresponding operators $\{B^k\}$ defined in \eqref{eq:B_k} satisfy
\[\rank(B^k) =\rank\bigl(\mathrm{Hess}_{\St}\Psi(\bV^*)\bigr)
\]
for all sufficiently large $k$.
\end{assumption}

\begin{proposition}\label{KKT}
Let $\{ \bV^k\}$ be a convergent sequence in $\St(K,R)$ with limit point $\bV^*$ such that $\nabla_{\St} \Psi(\bV^*) = 0$. Let $g^k \in \partial_{\St} \phi( \bV^k)$ be a collection of Riemannian subgradients such that $g^{k} \rightarrow g^* \in \partial_{\St} \phi(\bV^*)$ and  $\norm{g^{k} + \beta_{k} \nabla_{\St} \Psi( \bV^{k})} \rightarrow 0$ for some $ \beta \rightarrow \infty$. Suppose Assumption~\ref{ass:Psi2} and \ref{ass:rank} for $\Psi$ and $\{ \bV^k\}$ are satisfied. 
Then there exists $\omega^* \in \R^{R \times K}$ such that 
\[ g^* + \mathrm{Hess}_{\St} \Psi (\bV^*) \omega^* = 0.  \]
\end{proposition}

\begin{lemma}[Corollary 3.5 of Stewart~\cite{Stewart_1977}] \label{stewart}
Let $\{A_t\}$ be a sequence of matrices converging to $A$ as $t \to \infty$.
Then a necessary and sufficient condition for
\[
\lim_{t\to\infty} A_t^{+} = A^{+}
\]
is that $\operatorname{rank}(A_t) = \operatorname{rank}(A)$ for all sufficiently large $t$.
\end{lemma}

\begin{proof}[Proof of Proposition~\ref{KKT}]
For sufficiently large $k \in \N$, consider the tangent vectors $v_k = \mathrm{Log}_{\bV^*}(\bV^{k}) \in T_{\bV^*}\St$ and the operators $\{B^k\}$ as defined in \eqref{eq:B_k}. 
Defining 
\[
\omega_k := \beta_k v_k \in T_{\bV^*}\St ,
\]
we obtain by Taylor's theorem \eqref{eq:taylor} that $\beta_{k} \nabla_{\St} \Psi(\bV^{k})= P_{v_k} B^{k}\omega_k$. Hence, $g^{k}+P_{v_k} B^k\omega_k=g^{k}+\beta_{k}\nabla_{\St}\Psi(\bV^{k})$.
Next, define
\[\widehat{\omega}_k:=-(P_{v_k} B^k)^+ g^{k},\]
where $(\cdot)^+$ denotes the Moore--Penrose pseudoinverse. By definition of the pseudoinverse, $P_{v_k} B^k\widehat{\omega}^k$ is the orthogonal projection of $-g^{k}$ onto the range of $P_{v_k} B^k$. Consequently, 
\[\left\|g^{k}+P_{v_k} B^k \widehat{\omega}_k\right\|
\leq \left\|g^{k}+P_{v_k} B^k \omega_k\right\|
=
\left\|g^{k}+\beta_{k}\nabla_{\St}\Psi(\bV^{k})\right\|\to 0.\]
By the constant rank condition and Lemma \ref{stewart} we have that indeed $(P_{\bV^{k}} B^k )^+\to\big(\mathrm{Hess}_{\St}\Psi(\bV^*)\big)^+$. 
Combined with $g^{k} \to g^*$, this yields
\[\widehat{\omega}^k \to -\big(\mathrm{Hess}_{\St}\Psi(\bV^*)\big)^+ g^* \eqqcolon \omega^*.\]
Passing to the limit in $g^{k}+P_{\bV^{k}} B^k\widehat{\omega}^k\to 0$ finally gives
$g^* + \mathrm{Hess}_{\St}\Psi(\bV^*)\omega^* = 0$.
\end{proof}

\begin{theorem}
Let $\{\mathbf U^k\}$ be the sequence generated by \eqref{iterates} with arbitrary initialization.Set $\ga_k = \tilde \ga/(\beta_k \sqrt{k+1})$ with $\beta_k = C (k+1)^p$ and $p<1/4$. Define $\hat\UU^k = \prox_{\la \Omega_k/\beta_k}(\UU^k)$.
There is a subsequence $\{\tau(k)\}$ such that the following holds: 
\begin{enumerate}
    \item $U^{\tau(k)} \rightarrow \UU^*$
    \item $\hat U^{\tau(k)} \rightarrow \UU^*$
    \item $\UU^* \in N_\Psi$
    \item $ \exists \hat g^k \in \partial_{\St} \phi(\hat \UU^k) \colon \hat g^{\tau(k)} \rightarrow g^* \in \partial_{\St} \phi(\UU^*) \land \norm{\hat g^{\tau(k)} + \beta^{\tau(k)} \nabla_{\St} \Psi(\hat \UU^{\tau(k)})} \rightarrow 0.$
\end{enumerate}
Suppose that Assumptions~\ref{ass:Psi2} and~\ref{ass:rank} hold for $\Psi$ and the sequence $\{\hat\UU^{\tau(k)}\}$, respectively. If $N_\Psi$ is a smooth submanifold in a neighborhood of $\UU^*$, then $\UU^*$ is a stationary point of $\phi$ on $N_\Psi$, i.e.
\[0 \in \partial_{N_\Psi} \phi(\UU^*).\]
\end{theorem}

\begin{proof}
1. and 3. follow from Corollary \ref{cor:subseq}. 2. follows from Theorem \ref{thm:prox_conv} and 4. from Corollary \ref{cor:comb_conv} and Lemma \ref{conv_subgradient}.
Since $N_\Psi$ is a smooth submanifold in a neighborhood of $\UU^*$, for every $v \in T_{\UU^*}N_\Psi$ there exists a smooth curve
$c \colon (-\eps,\eps) \rightarrow N_\Psi$ such that
$c(0)=\UU^*$ and $c'(0)=v$. Since $\Psi$ is twice continuously differentiable and$\nabla_{\St}\Psi(c(t))=0$ for all $t$ in a neighborhood of zero, we have
\[\mathrm{Hess}_{\St}\Psi(\UU^*)v = \mathrm{Hess}_{\St}\Psi(\ga(0))\ga'(0) = \left. \frac{d}{dt}\nabla_{\St}\Psi(\ga(t))\right|_{t=0} = 0.\]
Hence, $T_{\UU^*}N_\Psi \subseteq \ker(\mathrm{Hess}_{\St} \Psi(\UU^*))$.
Applying Proposition  \ref{KKT} to the sequence $\{\hat \UU^{\tau(k)}\}$ gives the existence of $\omega^*$ such that $g^* + \mathrm{Hess}_{\St} \Psi (\UU^*) \omega^* = 0$. Consequently, 
\[g^* \in \range(\mathrm{Hess}_{\St} \Psi(\UU^*)) = \ker(\mathrm{Hess}_{\St} \Psi(\UU^*))^\perp \subset (T_{\UU^*}N_{\Psi})^\perp. \]
Therefore, $0 \in \partial_{N_\Psi} \phi(\UU^*)$.
\end{proof}

\section{Numerical Examples and Application}
\label{sec:Num}
In the numerical experiments, we investigate two aspects. First, we consider an artificial example to illustrate the algorithmic behavior, verify the derived rates, and empirically study the choice of hyperparameters $K$, $\bD$, and $\tilde\ga$. Second, we apply the algorithm to real data to obtain a diagonal frame decomposition of the Radon transform that promotes sparsity in the given data set. We then use this decomposition to recover images from noisy sinograms using regularizing filters.

\subsection{Visualization of theoretical results}

\begin{figure}
    \centering
    \begin{subfigure}[b]{0.48\textwidth}
        \centering
        \begin{subfigure}[t]{\linewidth}
            \centering
            \includegraphics[width=\linewidth]{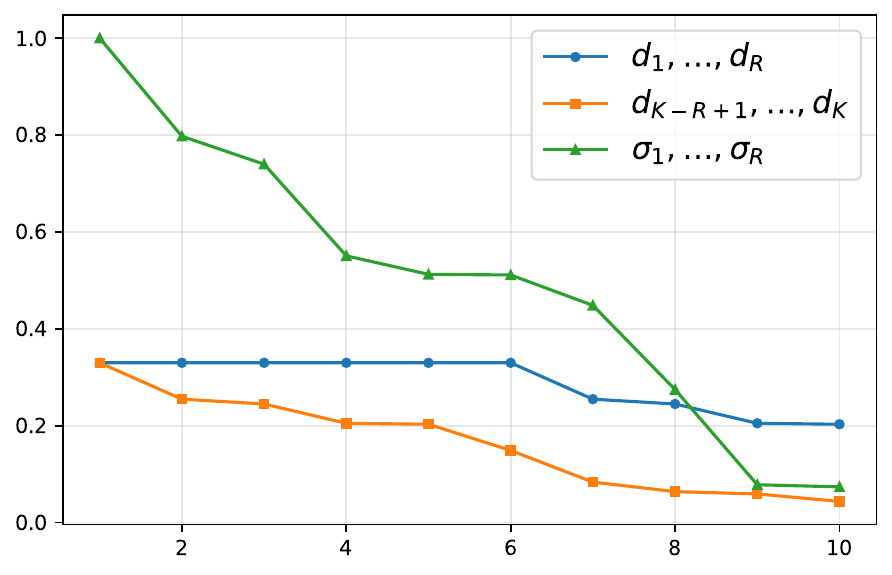}
            \caption{Heuristic choice of $\bD$.}
            \label{fig:d}
        \end{subfigure}
        \vspace{0.3cm}
        \begin{subfigure}[t]{\linewidth}
            \centering
            \includegraphics[width=\linewidth]{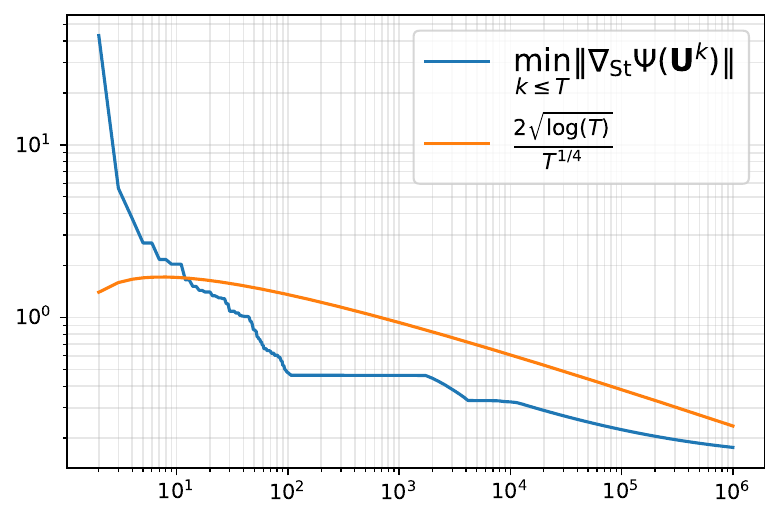}
            \caption{Convergence rate for the lower-level problem.}
            \label{fig:loss_grad}
        \end{subfigure}
    \end{subfigure}
    \hfill
    \begin{subfigure}[b]{0.48\textwidth}
        \centering
        \begin{subfigure}[t]{\linewidth}
            \centering
            \includegraphics[width=\linewidth]{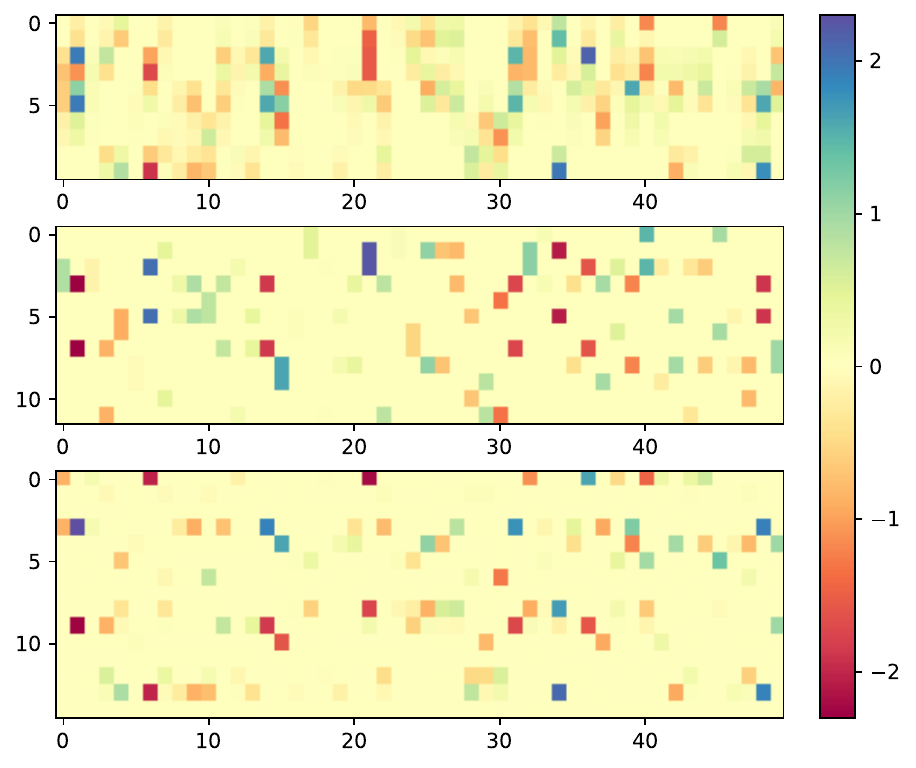}
            \caption{Top: data set $\bX$; middle: sparse representation with respect to the tight wavelet frame $\mathbf{W}$; bottom: sparse representation with respect to the algorithmic output $\UU^*$.}
            \label{fig:X}
        \end{subfigure}
        \vspace{0.3cm}
        \begin{subfigure}[t]{\linewidth}
            \centering
            \includegraphics[width=\linewidth]{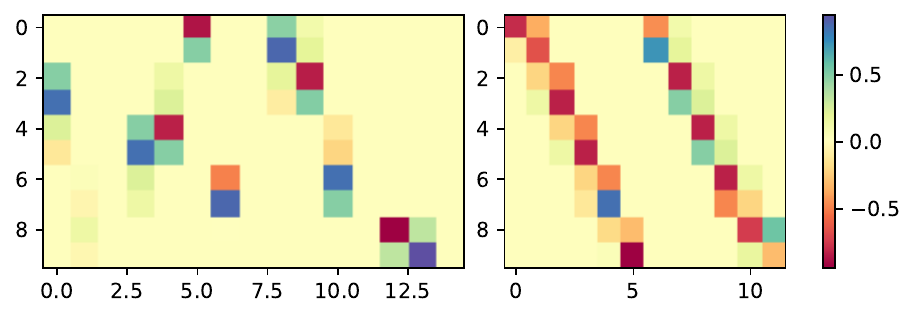}
            \caption{Left: algorithmic output $\UU^*$; right: tight wavelet frame $\mathbf{W}$.}
            \label{fig:U}
        \end{subfigure}
    \end{subfigure}
    \caption{Numerical results for the synthetic data set. The plots illustrate the choice of the diagonal matrix $\bD$, the convergence of the lower-level problem, the sparsity induced by the learned frame $\UU^*$, and a comparison of the learned frame with the wavelet frame $\mathbf{W}$ used to generate the sparse data.}
    \label{fig:U_new}
\end{figure}

In the first numerical experiment, we consider a randomly generated positive diagonal matrix $\Sigma \in \mathbb{R}^{10\times10}$ and seek a diagonal frame decomposition $(\UU,\bV,\bD)$ by first selecting $\bD$ and solving problem \eqref{pen_opt} via Algorithm~\ref{alg:RPS}. The data set $\bX$ consists of $S=50$ synthetically generated vectors in $\mathbb{R}^{10}$, each of which is $2$-sparse with respect to a tight wavelet frame $\mathbf{W}$ (i.e., the transformed data $\mathbf{W}^{\top}\bX$ has only two nonzero entries per column). The wavelet frame $\mathbf{W}$ is treated as unknown and is not used in the computations. The resulting data set and its transformation with respect to $\mathbf{W}$ is shown in Figure~\ref{fig:X} (top, middle). While $\|\bX\|_1\approx 118$, the sparsity measure with respect to the underlying wavelet frame is $\|\mathbf{W}^{\top}\bX\|_1\approx 85$.

We first set $K=15$ and choose $\bD$ according to the heuristic rule described in eq. \eqref{eq:d}. The choice of $\bD$ does not affect the algorithmic behavior, but determines the sparsity that can be achieved and the condition of $\bV$. Figure~\ref{fig:d} shows the diagonal entries of $\bD$ together with those of $\Sigma$. The ten largest and ten smallest values of $d$ relative to the corresponding singular values $\sigma$ are displayed to assess the applicability of the Fan--Pall theorem. Since the largest values of $d$ are smaller than the corresponding values of $\sigma$, the Fan--Pall conditions are violated, and hence no DFD exists for which both $\UU$ and $\bV$ are tight. Although $\scrC_{\bD}$ (defined in \eqref{eq:min-U-s.t.-C}) is empty, we can still approximate frames using Algorithm \ref{alg:RPS}, for which we know the performance guarantees worked out in the previous section. Throughout all experiments, we choose $\beta_k=(1+k)^{0.23}$, $\ga_k = \tilde \ga/(\beta_k \sqrt{k+1})$ for $\tilde \ga = 0.023$ unless stated otherwise. We initialize $\UU^0$ by applying the polar retraction to a random Gaussian matrix.

The output frame $\UU^*$ of the algorithm is shown in Figure~\ref{fig:U} (left), together with the wavelet frame $\mathbf{W}$ (right). Visually, the algorithm recovers the wavelet structures in which the data set is sparse. Moreover, the transformed data set $(\UU^*)^\top\bX$ is shown in Figure~\ref{fig:X} (bottom). It mostly exhibits the expected $2$-sparse structure, both visually and quantitatively, with $\|(\UU^*)^{\top}\bX\|_1 \approx 84$ which is even slightly smaller than the sparsity obtained with the original wavelet frame. 

\begin{figure}
    \centering
    \begin{subfigure}[b]{\textwidth}
        \centering
        \begin{subfigure}[t]{0.48\linewidth}
            \centering
            \includegraphics[width=\linewidth]{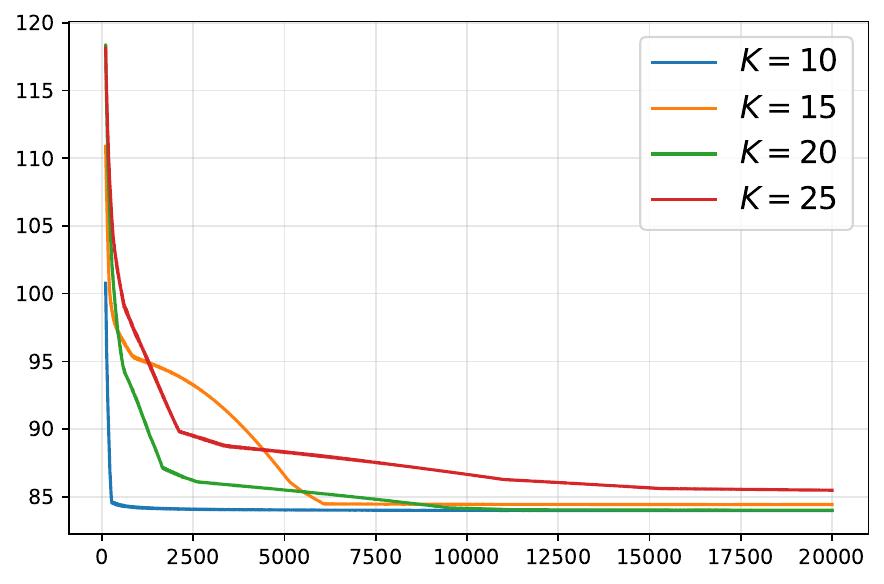}
        \end{subfigure}
        \hfill
        \begin{subfigure}[t]{0.48\linewidth}
            \centering
            \includegraphics[width=\linewidth]{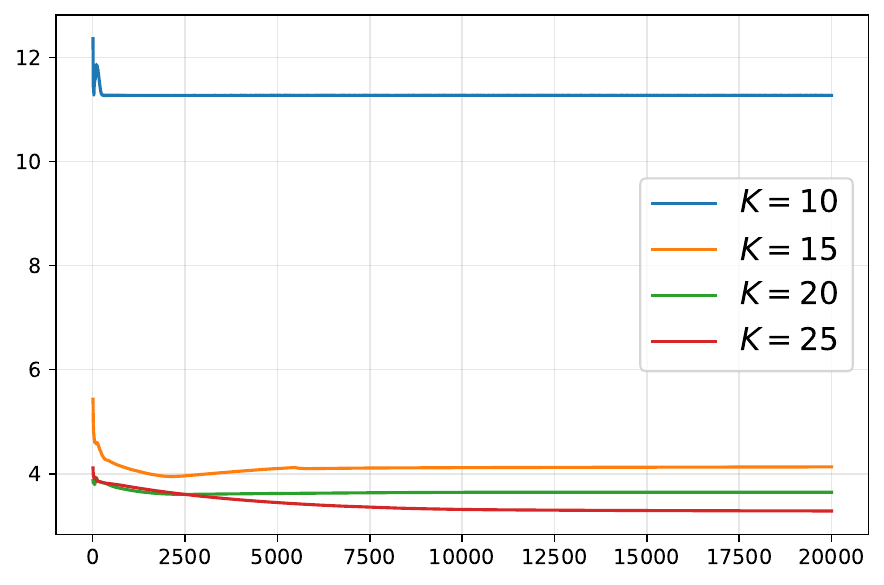}
        \end{subfigure}
    \end{subfigure}
    \caption{Sparsity loss $\norm{\UU^k \bX}_1$ (left) and condition number of $\bV^k = \Sigma^{-1}\UU^k\bD$ (right) for different choices of $K$ as a function of the iteration number for the synthetic example.}
    \label{fig:diff_K}
\end{figure}

As detailed in Subsection~\ref{sec:Pre_D}, the choice of $K$ is less critical than it may seem, provided that $K$ is sufficiently large, since some columns of $\UU$ may vanish. In our example, the wavelet frame that generates the sparse data has $12$ columns, so even for $K=15$, some columns of the learned frame are nearly zero. Figure~\ref{fig:diff_K} shows the loss $\norm{\UU^k \bX}_1$ (left) and the condition number of $\bV^k = \Sigma^{-1}\UU^k\bD$ (right) for different choices of $K$, with $\bD$ chosen according to \eqref{eq:d}.
The same level of sparsity is achieved for all choices of $K$, although $K=25$ requires more iterations or a larger step size. However, $K$ has a substantial impact on the condition number of $\bV^k$. Following the relaxation introduced in Section~\ref{sec:Opt}, the goal is therefore to obtain both a sparse representation $\UU^{\top}\bX$ and frame bounds of $\bV$ close to $1$. Since $\cond(\Sigma)\approx13.6$ in this example and $\bV$ cannot be tight for the chosen $\bD$, sufficiently large $K$ yields condition numbers around $3.6$ to $4$, whereas small $K$ results in condition numbers close to that of $\Sigma$, reducing the stability advantage of the DFD.

\begin{figure}
    \centering
    \begin{subfigure}[b]{\textwidth}
        \centering
        \begin{subfigure}[t]{0.48\linewidth}
            \centering
            \includegraphics[width=\linewidth]{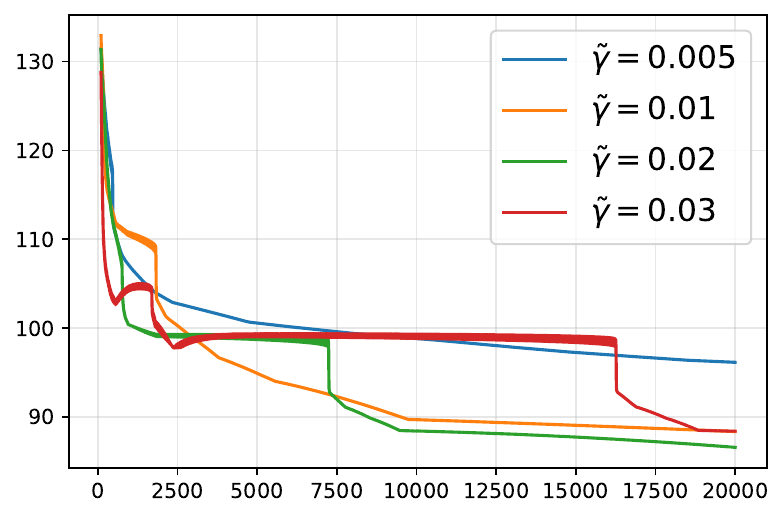}
        \end{subfigure}
        \hfill
        \begin{subfigure}[t]{0.48\linewidth}
            \centering
            \includegraphics[width=\linewidth]{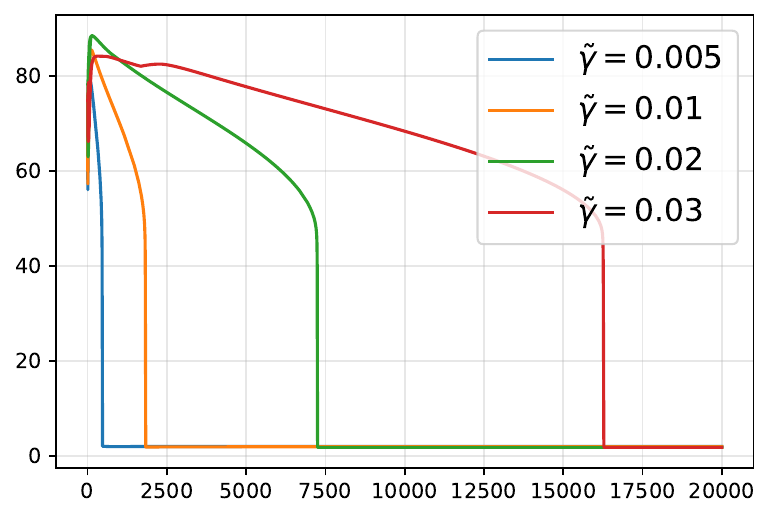}
        \end{subfigure}
    \end{subfigure}
    \caption{Sparsity loss $\norm{\UU^k \bX}_1$ (left) and condition number of $\bV^k = \Sigma^{-1}\UU^k\bD$ (right) for different step size constants $\tilde\ga$ for the synthetic example.}
    \label{fig:diff_gamma}
\end{figure}

Another observation concerns the sensitivity to the initial step size $\tilde\gamma$. In Figure \ref{fig:diff_gamma} we set $K=10$ and $\bD=\Sigma$ since the following observations are independent of the choice of $\bD$. Varying only the step size constant $\tilde \ga$ leads to substantial differences in the convergence of both the sparsity loss $\norm{\UU^k \bX}_1$  and the condition number of $\bV^k$. If the step size is too large, both the sparsity loss and the condition number remain high for a large number of iterations and then decrease rapidly as the algorithm approaches a good stationary point. If the step size is too small, convergence is excessively slow, and the desired sparsity may not be reached within a practical number of iterations. Thus, a suitable choice of $\tilde\gamma$ is crucial and also depends strongly on the condition numbers of $\Sigma$ and $\bD$.

\begin{figure}
    \centering
    \begin{subfigure}[b]{\textwidth}
        \centering
        \begin{subfigure}[t]{0.48\linewidth}
            \centering
            \includegraphics[width=\linewidth]{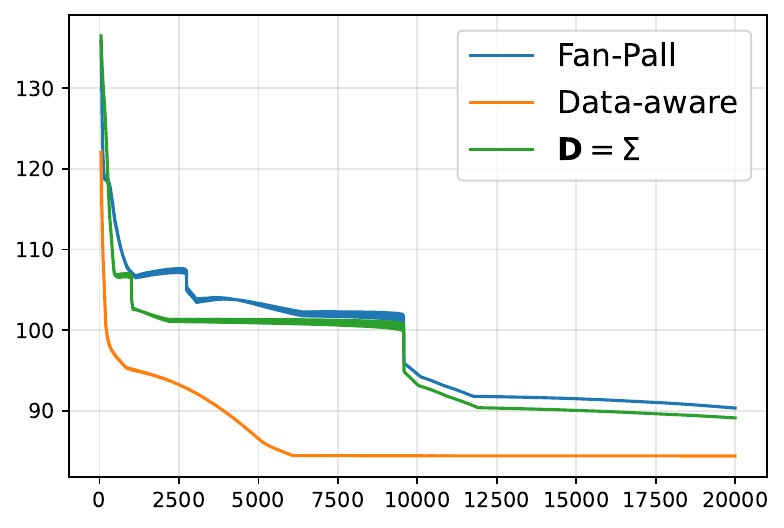}
        \end{subfigure}
        \hfill
        \begin{subfigure}[t]{0.48\linewidth}
            \centering
            \includegraphics[width=\linewidth]{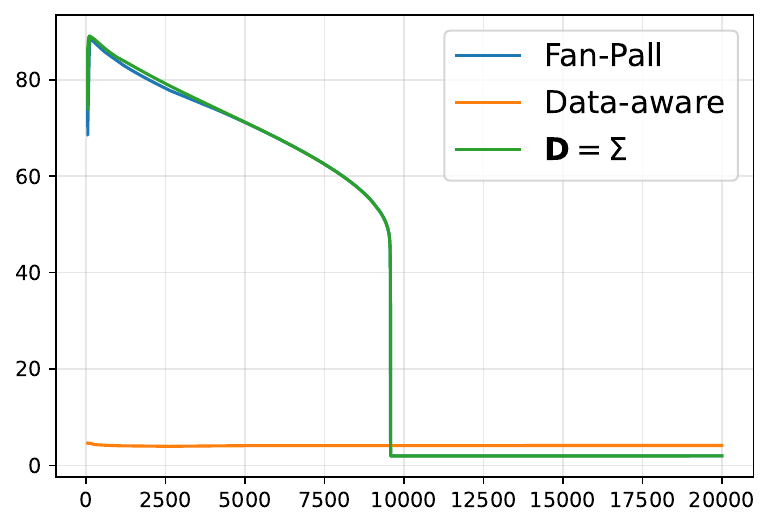}
        \end{subfigure}
    \end{subfigure}
    \caption{Sparsity loss $\norm{\UU^k \bX}_1$ (left) and condition number of $\bV^k = \Sigma^{-1}\UU^k\bD$ (right) for different choices of $\bD$ as a function of the iteration number for the synthetic example. In the Fan--Pall case $\bD$ is chosen according to \eqref{eq:D_fan_pall_short}, in the data-aware case according to \eqref{eq:d}, and in the third case we set $\bD=\Sigma$.}
    \label{fig:diff_D}
\end{figure}

Finally, the choice of $\bD$ has a clear influence on the resulting decomposition. In Figure~\ref{fig:diff_D}, we illustrate the convergence behavior for three choices: the Fan--Pall choice from \eqref{eq:D_fan_pall_short}, the data-aware heuristic choice from \eqref{eq:d}, and $\bD=\Sigma$. The data-aware choice achieves a lower sparsity loss than the Fan--Pall choice, but results in a larger condition number for $\bV^k$. In the Fan--Pall case, the condition number converges to approximately $2$, whereas in the data-aware case it stabilizes at around $4$. The choice $\bD=\Sigma$ exhibits a behavior similar to the Fan--Pall case. Although sparsity is an important objective, a large condition number can lead to unstable reconstruction methods, and thus limit the practical usefulness of the resulting decomposition. 
Note that the initial plateau in the Fan--Pall and $\bD=\Sigma$ cases is due to the choice of the initial step size. We use the same $\tilde\ga=0.023$ for all choices of $\bD$, while the data-aware choice results in smaller gradients and thus smaller effective updates in the algorithm.

\subsection{Filter-regularization with Learned DFD}

In the second numerical experiment, we consider a discretization of the two-dimensional Radon transform for the matrix $\bA$. The Radon transform maps a compactly supported square-integrable function to its line integrals over affine lines \cite{natterer2001mathematics}. It provides the mathematical foundation for computed tomography (CT) \cite{natterer01mathematical} and is widely used in imaging applications \cite{Kuchment2013, ZHANG2021104264}. Its discretized version $\bA$ maps a vectorized image $x \in \R^{m_1 \cdot m_2}$ to a vectorized observation (sinogram) $y \in \R^{n_1 \cdot n_2}$. The inversion of the Radon transform is ill-posed, meaning that small perturbations in the sinogram can lead to large reconstruction errors. This makes regularization essential for obtaining stable reconstructions from noisy data.

Learning a DFD for large-scale problems can be computationally demanding, as it requires the repeated computation of SVDs of high-dimensional matrices. However, this learning step only has to be performed once for a given operator and training data set. Once the DFD has been computed, reconstructing an image from a sinogram requires only two matrix multiplications and the evaluation of the corresponding filter, making reconstruction computationally inexpensive.
To keep the computational cost manageable, we conduct the numerical experiments using images of size \(20\times 20\), corresponding to 400 pixels, and sinograms of size \(20 \times 30\). Although this resolution is relatively low for practical image reconstruction, it provides a controlled setting for evaluating the reconstruction performance and error characteristics of the proposed method.

We first compute the SVD of $\bA$ and truncate the smallest singular values, whose contribution to the filtered reconstructions is negligible over the considered range of noise levels. We then learn a DFD for the slightly truncated operator $\bA_t = \bQ \Sigma \bW^\top$ with $\sigma_1 \geq \sigma_2 \geq \cdots \geq \eps > 0$. This truncation substantially reduces the computational cost of the learning algorithm and allows for larger step sizes, while having only a minor effect on the resulting filtered reconstructions. We choose $K=2R$, where $R$ is the rank of $\bA_t$, and we choose $\bD$ according to \eqref{eq:D_fan_pall_short} to satisfy the Fan-Pall Theorem \ref{thm:Fan_Pall}.  We initialize $\UU^0$ by applying the polar retraction to a random Gaussian matrix and we use $\beta_k=(1+k)^{0.23}$ and $\ga_k = \tilde \ga/(\beta_k \sqrt{k+1})$ with $\tilde \ga \approx \snorm{\UU^0}/\snorm{\beta_0 \nabla_{\St}\Psi(\UU^0) + \partial_{\St} \phi(\UU^0)}$.

To learn a tight DFD, we use $S=500$ real lung CT images as training data \citep{Soares2020}. Since we do not have a reliable a priori rule for choosing the initial step size or a suitable stopping criterion, we choose the initial step size such that the algorithm exhibits a rapid decay in the sparsity loss and the condition number of $\bV^k$ within a reasonable number of iterations ($k\leq 50000$). The influence of the initial step size on the convergence behavior and the resulting frame is similar to that observed in the artificial experiments shown in Figure~\ref{fig:diff_gamma}. Then, we selected the last iterate $\UU^\ast = \UU^k$ after $50000$ iterations.

Using the learned frames $\UU^\ast$ and  $\bV^\ast = \Sigma^{-1}\UU^\ast\bD$, we construct a DFD for $\bA_t$ according to Lemma~\ref{lem:DFD_SVD_equivalence}. We then apply the filter reconstruction \eqref{def:filter_reg} with the soft-thresholding filter $f_\al$ defined in \eqref{def:soft_thresholding} to reconstruct images from noisy sinograms $y^\delta$. More precisely, we compute
\begin{equation}
\mathcal{R}^{\mathrm{DFD}}_{\al}(y^\delta) = \bW^\top \UU^\ast f_\al(\bD^{+}, (\bV^\ast)^\top \bQ^\top y^\delta).
\end{equation}
The parameter $\alpha$ is chosen (unrealistically) to minimize the reconstruction error to the ground truth.

\begin{figure}
    \centering
    \begin{subfigure}[t]{0.32\textwidth}
        \centering
        \caption*{Relative $\ell^2$-error}
        \includegraphics[width=\linewidth]{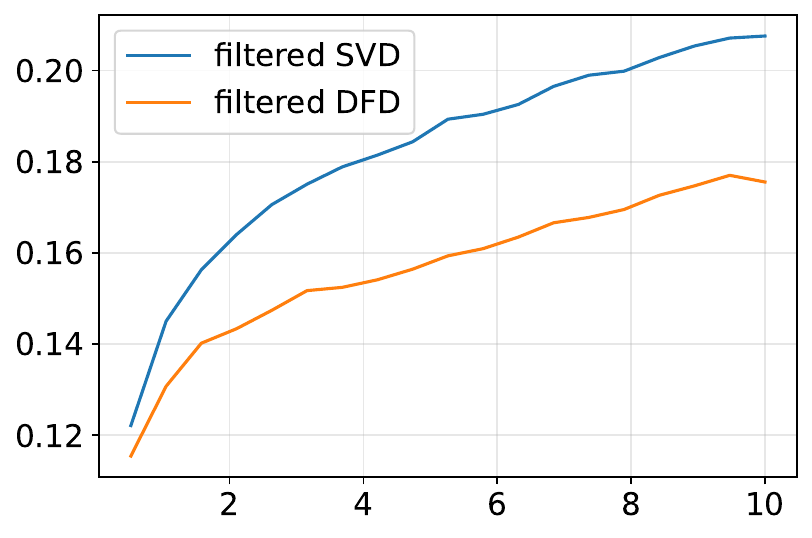}
    \end{subfigure}
    \hfill
    \begin{subfigure}[t]{0.32\textwidth}
        \centering
        \caption*{PSNR}
        \includegraphics[width=\linewidth]{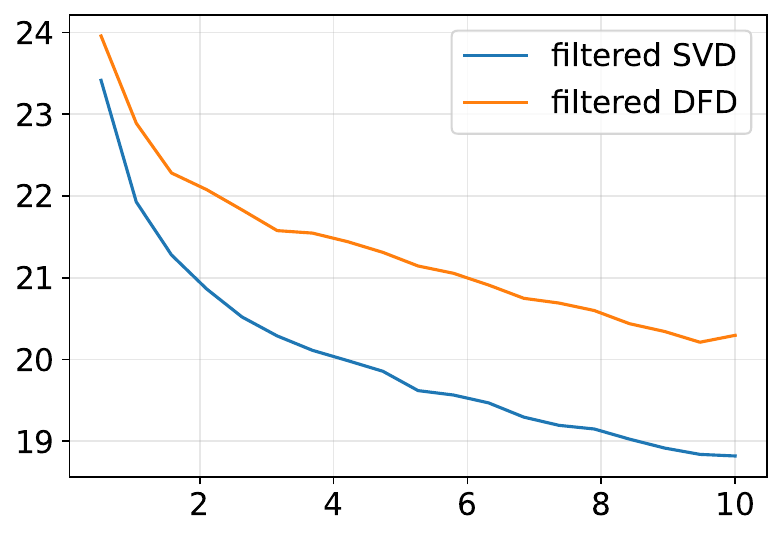}
    \end{subfigure}
    \hfill
    \begin{subfigure}[t]{0.32\textwidth}
        \centering
        \caption*{SSIM}
        \includegraphics[width=\linewidth]{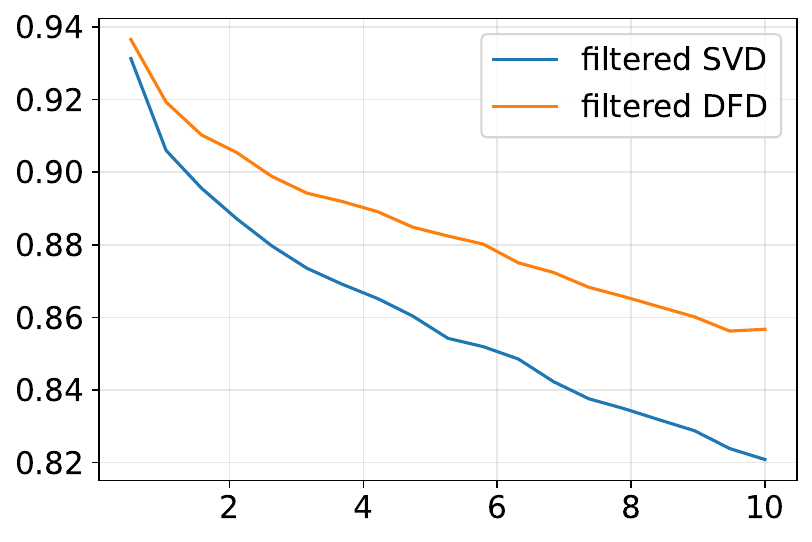}
    \end{subfigure}
    \caption{Comparison of the reconstruction quality measured by the relative $\ell^2$-error, PSNR, and SSIM. The $x$-axis displays the percentage of added Gaussian noise. For each noise level and regularization method, the parameter $\al$ is chosen to achieve the best possible reconstruction. The reported values are the mean over $50$ test images.}
    \label{fig:errors}
\end{figure}

To evaluate reconstruction quality, we use 50 independent test images and compute the mean relative $\ell^2$-error, peak signal-to-noise ratio (PSNR), and structural similarity index (SSIM) over the test set with respect to the ground truth. We generate noisy sinograms by adding Gaussian noise at different noise levels. For each noise level and reconstruction method, the regularization parameter $\al$ is chosen to yield the best reconstruction quality.
We compare the results with those obtained using the filtered SVD reconstruction
\begin{equation}
\mathcal{R}^{\mathrm{SVD}}_{\al}(y^\delta) = \bW^\top f_\al(\Sigma^{+}, \bQ^\top y^\delta). 
\end{equation}
In Figure \ref{fig:errors} we show the reconstruction quality achieved by filtered SVD and filtered DFD in different metrics. The filtered DFD with the data-aware learned frame consistently outperforms the filtered SVD for the considered noise levels and evaluation metrics, which demonstrates the benefit of adapting the frame to the training data. The performance of the learned DFD could potentially be further improved by increasing the number of training images and the number of iterations used in the learning algorithm.

\begin{figure}
    \centering
    \begin{subfigure}[t]{0.35\textwidth}
        \centering
        \includegraphics[width=\linewidth]{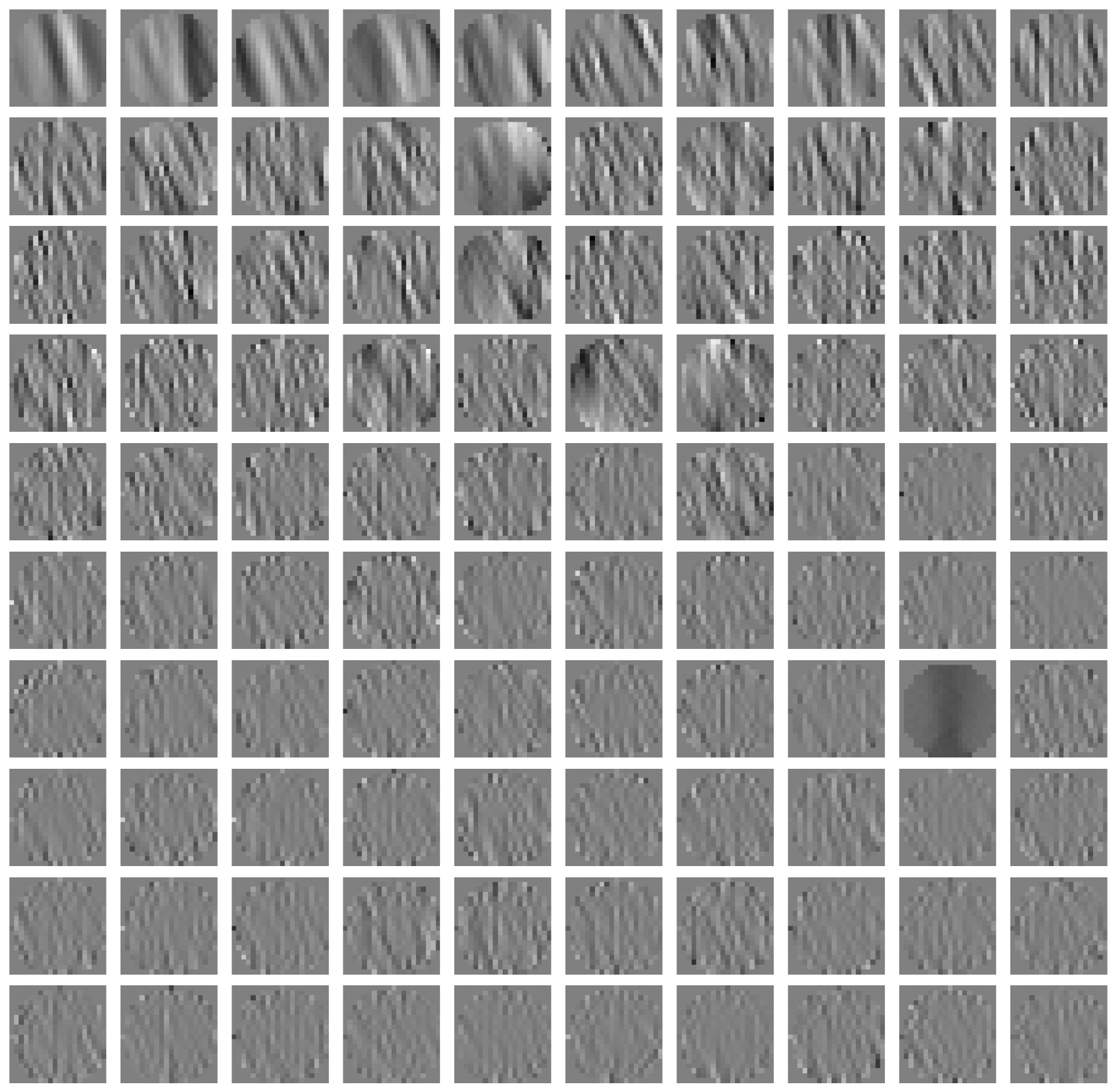}
    \end{subfigure}
   \hspace{0.5cm}
    \begin{subfigure}[t]{0.525\textwidth}
        \centering
        \includegraphics[width=\linewidth]{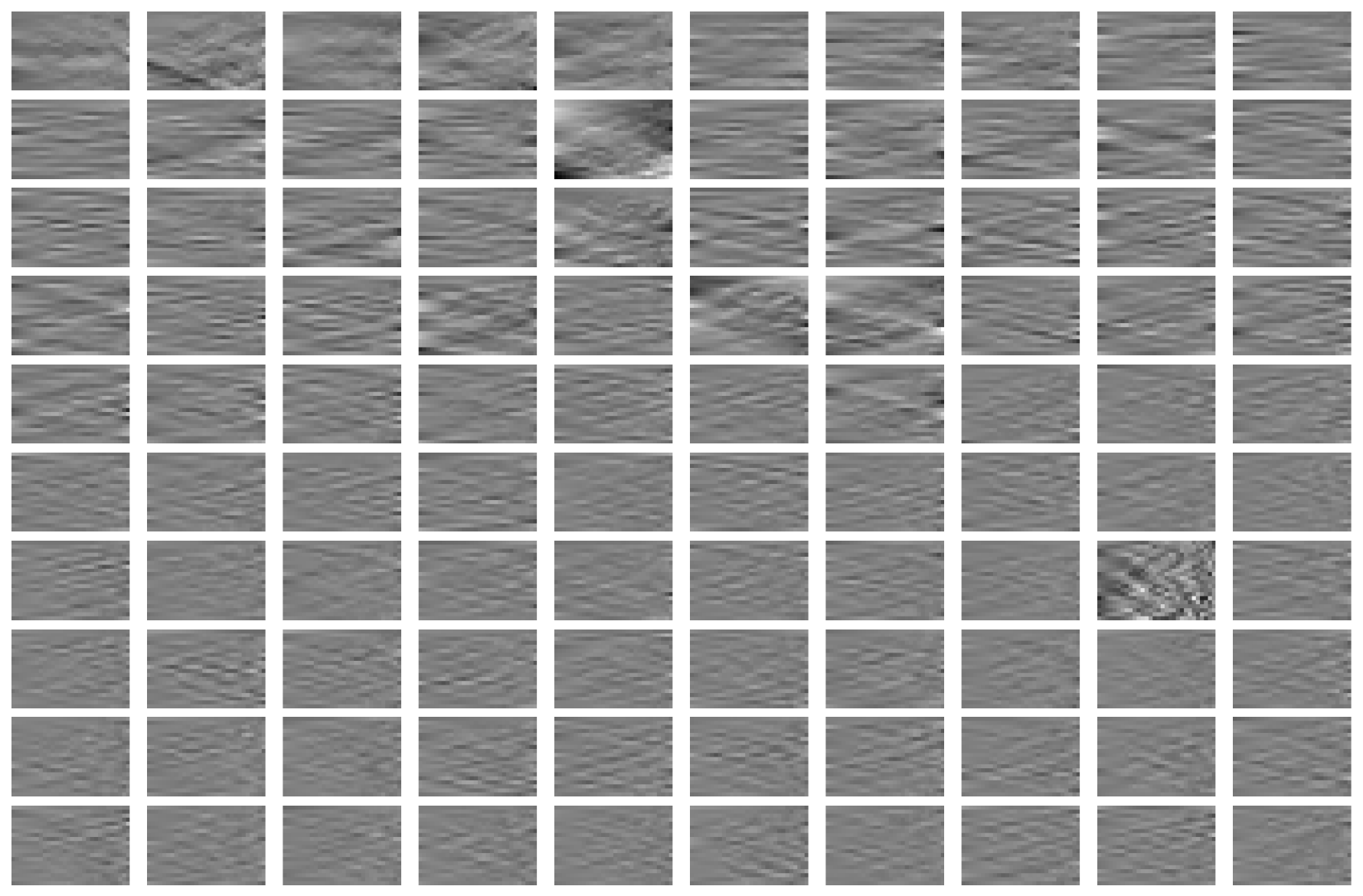}
    \end{subfigure}
    \caption{First 100 left DFD elements ($\UU$, left) and First 100 right DFD elements ($\bV$, right), arranged in a \(10\times10\) grid. Each element is reshaped into a \(20\times20\) image or a \(20\times 30\) sinogram, respectively. The grayscale is centered at zero, with negative values shown in darker shades and positive values in lighter shades.}
    \label{fig:frame_elements}
\end{figure}

To illustrate the structure of the learned DFD, we visualize the first 100 left and right frame elements in Figure \ref{fig:frame_elements}. The left frame elements are reshaped into \(20\times 20\) images, whereas the right frame elements are reshaped into \(20\times 30\) sinograms. All elements are displayed using the same grayscale scale, with zero corresponding to middle gray, negative values shown in darker shades, and positive values in lighter shades.

\section{Conclusion}
In this work, we introduced a data-aware approach for constructing diagonal frame decompositions whose frames in the signal space are tight and promote sparse representations of the given data. We formulated the construction as a bilevel optimization problem over the Stiefel manifold, where the nonconvex geometry of the feasible set leads to a nonconvex bilevel problem.

To address this problem, we proposed a penalized subgradient method and established convergence results in a general setting. Under suitable assumptions, we showed that the iterates converge to stationary points at both optimization levels and the derived convergence rates in terms of appropriate stationary measures. The numerical experiments confirm the theoretical behavior of the proposed method and demonstrate its practical ability to construct data-aware DFDs. In particular, the resulting decompositions provide effective regularization tools when applied to inverse problems, where they yield competitive reconstruction results.

Future work will include investigating other sparsity-promoting data terms and their practical performance in the construction of data-aware DFDs. From a theoretical perspective, an interesting direction is to extend convergence analysis towards near-stationary guaranties analogous to those established for nonconvex bilevel optimization in the Euclidean setting \cite{cortild2026nonconvexbilevel,NEURIPS2025_nonconvexbilevel}. On the computational side, parallelization and more efficient retractions on the Stiefel manifold could further improve the practical performance of the proposed method.

\section*{Declaration of generative AI and AI-assisted technologies in the manuscript preparation process}

During the preparation of this work, the authors used ChatGPT for language editing, brainstorming, and writing and debugging Python code. After using this tool, the authors reviewed and edited the content as needed and take full responsibility for the content of the published article.

\section*{Acknowledgements}
This research was funded in whole or in part by the Austrian Science Fund (FWF) J4954-N. 
For the purpose of open access, the author has applied a CC BY public copyright 
license to any Author Accepted Manuscript (AAM) version arising from this submission.

\appendix

\section*{Appendices}
\section{Tools on the Stiefel manifold}\label{app:Riemann}
In this appendix we collect some elements from manifold theory which are needed for our analysis. See \cite{Boumal_2023} for a general introduction to Riemannian optimization. 

\subsection{Basic notions}
\label{app:basic}
For $R\leq K$, we let $\R^{K\times R}$ be the Euclidean vectors space with inner product $\inner{\bA}{\bB}=\tr(\bA^{\top}\bB)$, and corresponding norm $\norm{\bA}_{F}:=\sqrt{\inner{\bA}{\bA}}$, the Frobenius norm. The Stiefel manifold on $\R^{K \times R}$ is  
\[
\St(K,R) \coloneqq  \{\mathbf{B} \in  \R^{K \times R} \mid \mathbf{B}^{\top} \mathbf{B} = \Id_{R} \}.
\]
We interpret this set as an embedded submanifold of $\R^{K\times R}$ via the locally defining function ${h:\R^{K\times R}\to\Sym(R), \bA\mapsto h(\bA)=\bA^{\top}\bA-\Id_{R}}$. Indeed, $h^{-1}(0)=\St(K,R)$ and $Dh(\bA)[\bV]=\bQ$ for $\bV=\frac{1}{2}\bA\bQ$ and $\bQ\in\Sym(R)$ is a surjective mapping. To shorten the notation we often omit the dimension and write $\St$. The tangent space of the Stiefel manifold at $\UU \in \St(K,R)$ is accordingly defined as 
\[
T_{\bA}\St(K,R)=\ker Dh(\bA)=\{\bV\in\R^{K\times R}\vert \bV^{\top}\bA+\bA^{\top}\bV=0\}. 
\]
This is an Euclidean vector space inheriting the scalar product of its ambient space $\R^{K\times R}$. 

For a general embedded submanifold $\scrM$, the product set $T\scrM=\{(x,\xi)\vert x\in\scrM,\xi\in T_{x}\scrM\}$ is its \emph{tangent bundle}. 

\subsection{Retractions}
A \emph{retraction} on a manifold $\scrM$ is a smooth map 
\[
\Retract:T\scrM\to\scrM,\quad (\bA,\xi)\mapsto \Retract_{\bA}(\xi)
\]
such that each curve $c(t)=\Retract_{\bA}(t\xi)$ satisfies $c(0)=\bA$ and $c'(0)=\xi$. There are multiple ways to define retractions on a given manifold, and the practical choice is usually guided by the geometric and algebraic structure of $\scrM$. For the Stiefel manifold $\scrM=\St(K,R)$ various efficiently computable retractions are known. A common choice is the $Q$-factor retraction: $\Retract_{\bA}(\xi)=\bQ$, where $\bQ\mathbf{R}=X+\xi$ is a thin QR decomposition.
A particular efficient retraction, also used in this paper, is the polar retraction 
\[
\Retract_{\bA}(\xi)=(\bA+\xi)^{\top}(\Id_{R}+\xi\xi^{\top})^{-1/2}.
\]
It can be efficiently computed through the SVD of $\bA+\xi$. To wit, if $\bA+\xi=\bZ\Lambda\bW^{\top}$, then $\Retract_{\bA}(\xi)=\bZ\bW^{\top}$. 
An important property of the polar retraction is its nonexpansiveness.
\begin{lemma}[Lemma 1, \cite{li2021weakly}]
\label{polar_retraction}
Let $\UU \in \St(K,R)$ and let $\xi \in T_{\UU}\St$. Define $\UU_+ \coloneqq \UU + \xi$. Then the polar decomposition based retraction~(2.4) satisfies
\[
\Retract_{\UU}(\xi)
=
\UU_+
\bigl(\UU_+^\top \UU_+\bigr)^{-1/2}
=
\proj_{\St}(\UU_+).
\]
Moreover,
\begin{equation}\label{eq:Retraction}
\|\operatorname{Retr}_{\UU}(\xi)-\bar{\UU}\|_F
\le
\|\UU_+-\bar{\UU}\|_F\qquad \forall \bar{\UU}\in\St(K,R).
\end{equation}
\end{lemma}

\subsection{Derivatives}
For a smooth map $F\colon\scrM\to\scrM'$ btween manifolds, its differential at $\bA\in\scrM$ is the linear map
$\mathrm{D}F(\bA)\colon T_{\bA}\scrM\to T_{F(\bA)}\scrM'$ defined by
\[
\mathrm{D}F(\bA)[\bB]=\left.\frac{\mathrm{d}}{\mathrm{d}t}F(c(t))\right|_{t=0},
\]
where $c$ is any smooth curve on $\scrM$ satisfying $c(0)=\bA$ and $c'(0)=\bB$; see \cite[Def.~3.34]{Boumal_2023}. For a smooth function $f\colon\scrM\to\R$, the Riemannian gradient $\grad_{\scrM}f(\bA)\in T_{\bA}\scrM$ is the unique tangent vector satisfying $\mathrm{D}f(\bA)[\bB]=\langle \bB,\grad_{\scrM}f(\bA)\rangle_{\bA}$ for all $\bB\in T_{\bA}\scrM$;
see \cite[Def.~3.58]{Boumal_2023}. The Riemannian Hessian is the covariant derivative of the Riemannian gradient and can be viewed as a linear map
\[
\mathrm{Hess}_{\scrM}f(\bA)\colon T_{\bA}\scrM\to T_{\bA}\scrM.
\]

Since the Stiefel manifold is an embedded submanifold of the Euclidean space $\R^{K\times R}$, Riemannian derivatives can be expressed in terms of Euclidean derivatives and the orthogonal projection onto the tangent space. 
The orthogonal projection onto the tangent space $T_{\bA}\St(K,R)$ is the operator such that $\bV-\proj_{T_{\bA}\St}(\bV)\in T_{\bA}\St(K,R)^{\bot}$. Using the Euclidean structure of the embedding space, one can easily compute that 
\[
\proj_{T_{\bA}\St}(\bV)=\bV-\bA\frac{\bA^{\top}\bV+\bV^{\top}\bA}{2}=(\Id_{R}-\bA\bA^{\top})\bV+\bA\frac{\bA^{\top}\bV-\bV^{\top}\bA}{2}.
\]
The projector yields a convenient formula for the  Riemannian gradient of a differentiable function $f \colon \St \rightarrow \R$ in terms of a differentiable extension $\bar{f} \colon \R^{K\times R}\to\R$ by \cite[][Proposition 3.61]{Boumal_2023}
\[
\nabla_{\St}f(\bA)=\proj_{T_{\bA}\St}\bigl(\nabla\bar{f}(\bA)\bigr)=\nabla \bar{f}(\bA)-\bA\frac{\bA^{\top}\nabla\bar{f}(\bA)+\nabla\bar{f}(\bA)^{\top}\bA}{2}.
\]
If $f$ is smooth enough the Riemannian Hessian $\mathrm{Hess}_{\St} f(\bA)$ is a self-adjoint linear map on $T_{\bA}\St$ given by  
\[ 
\mathrm{Hess}_{\St} f(\bA) [\bB] = \proj_{T_{\bA}\St}(\mathrm{Hess}_{\St} \bar f (\bA) [\bB]) + \mathcal{P}_{\bB} ( \proj_{T_{\bA}\St}^\perp( \nabla \bar f (\bA)))
\]
where $\mathcal{P}_{\bB}$ is the differential of $\bA \mapsto \proj_{T_{\bA}\St}$ at $\bA$ along $\bB$ and $\proj_{T_{\bA}\St}^\perp = \Id - \proj_{T_{\bA}\St}$. Since we use only first-order methods, the Hessian serves only as a theoretical tool, and we do not need to compute it explicitly in our experiments.

For a convex function $\bar g \colon \R^{K\times R}\to\R$, the Riemannian (Clarke) subdifferential $\partial_{\St}g(\bA)$ of the restriction $g = \bar g|_{\St}$ is also given by the projector \cite[Theorem 5.1]{yang2014optimality}
\[
\partial_{\St}g(\bA)=\proj_{T_{\bA}\St}(\partial \bar g (\bA))\qquad \forall \bA\in\St(K,R).
\]
In particular, although $g$ is generally no longer convex, any Euclidean subgradient $g'(\bA)\in\partial \bar g(\bA)$ gives rise to a Riemannian subgradient through the projection operator $\proj_{T_{\bA}\St}(g'(\bA))$.

\subsection{Parallel Transport and Taylor expansion}
\label{app:taylor}
In our optimality analysis, we make use of the Taylor expansion on the Stiefel manifold. We briefly introduce the basic tools required to apply it; for a detailed treatment, we refer to \cite[Chap.~5, Chap.~10]{Boumal_2023}. Notably, this tool is not specific to the Stiefel manifold but applies to arbitrary Riemannian submanifolds.

The so-called parallel transport is a tool that allows us to transport a tangent vector $\xi \in T_{\bA}\St$ from $\bA$ to another point $\bB \in \St$. Given a smooth curve $c \colon I \rightarrow \St$ with $c(t_0) = \bA$ and $c(t_1)=\bB$, there exists a unique smooth map (vector field on $c$) $Z \colon I \rightarrow T\St$ such that $Z(t) \in T_{c(t)}\St$ for all $t \in I$, $Z(t_0) = \xi$, and $\proj_{c(t)}(\frac{\dif}{\dif t}Z(t)) = 0$. Here, $\frac{\dif}{\dif t}Z$ denotes the classical derivative of $Z$ in the embedding space. A vector field $Z$ satisfying these requirements is said to be parallel.
The parallel transport of tangent vectors at $c(t_0)$ to the tangent space at $c(t_1)$ along $c$ is the map 
\[
PT^c_{t_1 \leftarrow t_0} \colon T_{c(t_0)} \St \rightarrow T_{c(t_1)}\St
\]
defined by $PT^c_{t_1 \leftarrow t_0}(\xi) = Z(t_1)$. In the case of the Stiefel manifold, parallel transport is a linear isometry and $PT^c_{t \leftarrow t}$ is the identity.

To obtain a second-order Taylor expansion, a particular curve between points is used, namely the geodesic generated by the so-called exponential map.
For every $(\bA, \xi) \in T \St$ there exists a unique maximal geodesic $\ga \colon I \rightarrow \St$ with $\ga(0) = \bA$ and $\ga'(0) = \xi$, where maximal means that the interval $I$ is as large as possible. If $I$ contains $[0,1]$ the exponential map is defined by 
\[
\mathrm{Exp}(\bA,\xi) = \mathrm{Exp}_\bA(\xi) = \ga(1).
\]
For $\bB$ sufficiently close to $\bA$, there exists a tangent vector
$\xi\in T_{\bA}\St$ such that $\bB=\mathrm{Exp}_{\bA}(\xi)$. 
This tangent vector is denoted by the logarithmic map
\[
    \mathrm{Log}_{\bA}(\bB)=\xi.
\]
Now consider the geodesic $c(t) = \mathrm{Exp}_{\bA}(t \xi)$ between $\bA$ and $\bB$. We write $P_{t \xi} = PT^c_{t \leftarrow 0}$. Then, if $f$ is at least twice continuously differentiable, the following Taylor expansion of the Riemannian gradient holds:
\[
P^{-1}_{t\xi} \nabla_{\St} f(c(t)) = \nabla_{\St} f(\bA) + \int_{0}^{t} ( P_{\tau \xi}^{-1} \circ \mathrm{Hess}_{\St} f(c(\tau)) \circ P_{\tau \xi})[\xi] \dif \tau.
\]

\section{Weakly convex functions}\label{app:weakly_convex}
A function $f:\mathcal{U} \subset \mathbb{R}^n\to\mathbb{R}$ is called $\tau$-weakly convex if
\[
x\mapsto f(x)+\frac{\tau}{2}\norm{x}^2
\]
is convex for some $\tau\geq 0$. Thus, weak convexity allows for a controlled deviation from convexity and includes convex functions as the special case $\tau=0$.
If $f$ is differentiable and $\tau$-weakly convex it satisfies the gradient inequality
\begin{equation}\label{eq:weakly}
f(\hat{x})\geq f(x)+\inner{\nabla f(x)}{\hat{x}-x}-\frac{\tau}{2}\norm{\hat{x}-x}^{2}\qquad\forall \hat{x},x\in\mathcal{U}. 
\end{equation}

\bibliographystyle{plain}
\bibliography{sample}

@article{Soares2020,
  author  = {Soares, Eduardo and Angelov, Plamen and Biaso, Sarah and
             Froes, Michele Higa and Abe, Daniel Kanda},
  title   = {{SARS-CoV-2 CT-scan dataset}: A large dataset of real patients CT scans for {SARS-CoV-2} identification},
  journal = {medRxiv},
  year    = {2020},
  doi     = {10.1101/2020.04.24.20078584}
}

@book{natterer2001mathematics,
    title={The mathematics of computerized tomography},
    author={Natterer, F.},
    year={2001},
    publisher={SIAM}
  }

@Book{natterer01mathematical,
  title     = {Mathematical Methods in Image Reconstruction},
  publisher = {SIAM},
  year      = {2001},
  author    = {Natterer, F. and W{\"u}bbeling, F.},
  volume    = {5},
  series    = {Monographs on Mathematical Modeling and Computation},
  address   = {Philadelphia, PA},
  }

@book{Kuchment2013,
  author    = {Kuchment, Peter},
  title     = {The Radon Transform and Medical Imaging},
  publisher = {SIAM},
  year      = {2013},
  doi       = {10.1137/1.9781611973297}
}

@article{ZHANG2021104264,
title = {A local radon transform for seismic random noise attenuation},
journal = {Journal of Applied Geophysics},
volume = {186},
pages = {104264},
year = {2021},
issn = {0926-9851},
doi = {https://doi.org/10.1016/j.jappgeo.2021.104264},
url = {https://www.sciencedirect.com/science/article/pii/S0926985121000112},
author = {Quan Zhang and Hang Wang and Wei Chen and Guangtan Huang},
}

@book{Ch02,
author = {Christensen, O.},
publisher = {Birkh{\"{a}}user},
title = {An Introduction to Frames and Riesz Bases},
year = {2002}
}

@article{TroppIEEE05,
	author = {Tropp, Joel A and Dhillon, Inderjit S and Heath, Robert W and Strohmer, Thomas},
	date = {2005},
	journal = {IEEE Transactions on information theory},
	number = {1},
	pages = {188-209},
	publisher = {IEEE},
	title = {Designing structured tight frames via an alternating projection method},
	volume = {51},
	year = {2005}}

@ARTICLE{Chaux2006,
  author={Chaux, C. and Duval, L. and Pesquet, J.-C.},
  journal={IEEE Transactions on Image Processing}, 
  title={Image analysis using a dual-tree M-band wavelet transform}, 
  year={2006},
  volume={15},
  number={8},
  pages={2397-2412},
  doi={10.1109/TIP.2006.875178}}

@ARTICLE{Selesnick05,
  author={Selesnick, I.W. and Baraniuk, R.G. and Kingsbury, N.C.},
  journal={IEEE Signal Processing Magazine}, 
  title={The dual-tree complex wavelet transform}, 
  year={2005},
  volume={22},
  number={6},
  pages={123-151},
  doi={10.1109/MSP.2005.1550194}}

@article{eldar2003sampling,
  title={Sampling with arbitrary sampling and reconstruction spaces and oblique dual frame vectors},
  author={Eldar, Yonina C},
  journal={Journal of Fourier Analysis and Applications},
  volume={9},
  number={1},
  pages={77--96},
  year={2003},
  publisher={Springer}
}

@article{blum2010sobolev,
  title={Sobolev duals in frame theory and sigma-delta quantization},
  author={Blum, James and Lammers, Mark and Powell, Alexander M and Y{\i}lmaz, {\"O}zg{\"u}r},
  journal={Journal of Fourier Analysis and Applications},
  volume={16},
  number={3},
  pages={365--381},
  year={2010},
  publisher={Springer}
}

@article{bodmann2007frame,
  title={Frame paths and error bounds for sigma--delta quantization},
  author={Bodmann, Bernhard G and Paulsen, Vern I},
  journal={Applied and Computational Harmonic Analysis},
  volume={22},
  number={2},
  pages={176--197},
  year={2007},
  publisher={Elsevier}
}

@article{candes2004new,
  title={New tight frames of curvelets and optimal representations of objects with piecewise C2 singularities},
  author={Cand{\`e}s, Emmanuel J and Donoho, David L},
  journal={Communications on Pure and Applied Mathematics: A Journal Issued by the Courant Institute of Mathematical Sciences},
  volume={57},
  number={2},
  pages={219--266},
  year={2004},
  publisher={Wiley Online Library}
}

@article{kutyniok2011compactly,
  title={Compactly supported shearlets are optimally sparse},
  author={Kutyniok, Gitta and Lim, Wang-Q},
  journal={Journal of Approximation Theory},
  volume={163},
  number={11},
  pages={1564--1589},
  year={2011},
  publisher={Elsevier}
}

@article{strohmer2003grassmannian,
  title={Grassmannian frames with applications to coding and communication},
  author={Strohmer, Thomas and Heath Jr, Robert W},
  journal={Applied and computational harmonic analysis},
  volume={14},
  number={3},
  pages={257--275},
  year={2003},
  publisher={Elsevier}
}

@book{daubechies1992ten,
  title={Ten lectures on wavelets},
  author={Daubechies, Ingrid},
  year={1992},
  publisher={SIAM}
}

@article{KUTYNIOK20132225,
title = {Scalable frames},
journal = {Linear Algebra and its Applications},
volume = {438},
number = {5},
pages = {2225-2238},
year = {2013},
issn = {0024-3795},
doi = {https://doi.org/10.1016/j.laa.2012.10.046},
url = {https://www.sciencedirect.com/science/article/pii/S0024379512007768},
author = {Gitta Kutyniok and Kasso A. Okoudjou and Friedrich Philipp and Elizabeth K. Tuley}
}

@book{Waldron2018,
  author       = {Shayne F. D. Waldron},
  title        = {An Introduction to Finite Tight Frames},
  series       = {Applied and Numerical Harmonic Analysis},
  year         = {2018},
  edition      = {1},
  publisher    = {Birkh{\"a}user},
  address      = {New York, NY},
  doi          = {10.1007/978-0-8176-4815-2},
  isbn         = {978-0-8176-4814-5},
  isbn2        = {978-1-4939-7920-2},
  eisbn        = {978-0-8176-4815-2},
  issn         = {2296-5009},
  eissn        = {2296-5017},
  pages        = {XX + 587},
  note         = {25 b/w illustrations, 12 colour illustrations},
  url          = {https://doi.org/10.1007/978-0-8176-4815-2}
}

@article{candes2002recovering,
  title={Recovering edges in ill-posed inverse problems: Optimality of curvelet frames},
  author={Candes, Emmanuel J and Donoho, David L},
  journal={The Annals of Statistics},
  volume={30},
  number={3},
  pages={784--842},
  year={2002},
  publisher={Institute of Mathematical Statistics}
}

@article{ebner2023regularization,
  title={Regularization of inverse problems by filtered diagonal frame decomposition},
  author={Ebner, Andrea and Frikel, J{\"u}rgen and Lorenz, Dirk and Schwab, Johannes and Haltmeier, Markus},
  journal={Applied and Computational Harmonic Analysis},
  volume={62},
  pages={66--83},
  year={2023},
  publisher={Elsevier}
}

@article{frikel2018efficient,
  title={Efficient regularization with wavelet sparsity constraints in photoacoustic tomography},
  author={Frikel, J{\"u}rgen and Haltmeier, Markus},
  journal={Inverse Problems},
  volume={34},
  number={2},
  year={2018},
  publisher={IOP Publishing}
}

@article{Fr19,
author = {Frikel, J{\"{u}}rgen and Haltmeier, Markus},
issn = {2297024X},
journal = {Trends in Mathematics},
pages = {163--178},
title = {{Sparse regularization of inverse problems by operator-adapted frame thresholding}},
year = {2020}
}

@article{ebner2023convergence,
year = {2024},
publisher = {IOP Publishing},
volume = {40},
number = {5},
author = {Andrea Ebner and Markus Haltmeier},
title = {Convergence of non-linear diagonal frame filtering for regularizing inverse problems},
journal = {Inverse Problems}
}

@article{donoho1995nonlinear,
  title={Nonlinear solution of linear inverse problems by wavelet--vaguelette decomposition},
  author={Donoho, David L},
  journal={Applied and computational harmonic analysis},
  volume={2},
  number={2},
  pages={101--126},
  year={1995},
  publisher={Elsevier}
}

@book{En96,
author = {Engl, Heinz W. and Hanke, Martin and Neubauer, Andreas},
booktitle = {Regularization of Inverse Problems},
isbn = {978-0-7923-4157-4},
publisher = {Springer},
title = {{Regularization of Inverse Problems}},
year = {1996}
}

@inproceedings{NEURIPS2025_nonconvexbilevel,
 author = {Cao, Jincheng and Jiang, Ruichen and Yazdandoost Hamedani, Erfan and Mokhtari, Aryan},
 booktitle = {Advances in Neural Information Processing Systems},
 doi = {10.52202/085713-1079},
 editor = {D. Belgrave and C. Zhang and H. Lin and R. Pascanu and P. Koniusz and M. Ghassemi and N. Chen},
 pages = {32126--32153},
 publisher = {Curran Associates, Inc.},
 title = {On the Complexity of Finding Stationary Points in Nonconvex Simple Bilevel Optimization},
 url = {https://proceedings.neurips.cc/paper_files/paper/2025/file/2e2ac035e8f2c4e357e6242805540423-Paper-Conference.pdf},
 volume = {38, Main Conference},
 year = {2025}
}

@misc{cortild2026nonconvexbilevel,
      title={Stochastic Nonconvex Bilevel Optimization: Improved Rates Without Rare-Visit Assumption}, 
      author={Daniel Cortild and Mathias Staudigl and Juan Peypouquet and Coralia Cartis},
      year={2026},
      eprint={2609.06580},
      archivePrefix={arXiv},
      primaryClass={math.OC},
      url={https://arxiv.org/abs/2609.06580}, 
}

@article{dempe2021simple,
  title={Simple bilevel programming and extensions},
  author={Dempe, Stephan and Dinh, Nguyen and Dutta, Joydeep and Pandit, Tanushree},
  journal={Mathematical Programming},
  volume={188},
  number={1},
  pages={227--253},
  year={2021},
  publisher={Springer}
}

@article{solodov2007explicit,
  title={An explicit descent method for bilevel convex optimization},
  author={Solodov, Mikhail},
  journal={Journal of Convex Analysis},
  volume={14},
  number={2},
  pages={227},
  year={2007},
  publisher={HELDERMANN VERLAG LANGER GRABEN 17, 32657 LEMGO, GERMANY}
}

@article{DavDrusSIOPT19,
author = {Davis, Damek and Drusvyatskiy, Dmitriy},
title = {Stochastic Model-Based Minimization of Weakly Convex Functions},
journal = {SIAM Journal on Optimization},
volume = {29},
number = {1},
pages = {207-239},
year = {2019},
doi = {10.1137/18M1178244},
URL = {https://doi.org/10.1137/18M1178244},
eprint = {https://doi.org/10.1137/18M1178244}
}

@article{10530399,
	author = {Staudigl, Mathias and Qu, Siqi},
	doi = {10.1109/LCSYS.2024.3401019},
	journal = {IEEE Control Systems Letters},
	pages = {622-627},
	title = {Tikhonov Regularized Exterior Penalty Dynamics for Constrained Variational Inequalities},
	volume = {8},
	year = {2024}}

@article{AttCzarPey11,
	author = {Attouch, H\'{e}dy and Czarnecki, Marc-Olivier and Peypouquet, Juan},
	doi = {10.1137/110820300},
	eprint = {https://doi.org/10.1137/110820300},
	journal = {SIAM Journal on Optimization},
	number = {4},
	pages = {1251-1274},
	title = {Coupling Forward-Backward with Penalty Schemes and Parallel Splitting for Constrained Variational Inequalities},
	url = {https://doi.org/10.1137/110820300},
	volume = {21},
	year = {2011}}

@article{PeyStaQu26,
	author = {Peypouquet, Juan and Qu, Siqi and Staudigl, Mathias},
	journal = {Applied Mathematics \& Optimization},
	number = {2},
	pages = {36},
	title = {Asymptotic Behavior of Penalty Dynamics for Constrained Variational Inequalities},
	volume = {93},
	year = {2026}}

@article{Att96,
	author = {Attouch, Hedy},
	journal = {SIAM Journal on Optimization},
	number = {3},
	pages = {769--806},
	publisher = {SIAM},
	title = {Viscosity solutions of minimization problems},
	volume = {6},
	year = {1996}}

@article{BaiCom01,
	author = {Baillon, Jean-Bernard and Cominetti, Roberto},
	journal = {Journal of Functional Analysis},
	number = {2},
	pages = {263--273},
	publisher = {Elsevier},
	title = {A convergence result for nonautonomous subgradient evolution equations and its application to the steepest descent exponential penalty trajectory in linear programming},
	volume = {187},
	year = {2001}}

@article{ChenSIOPT2020,
author = {Chen, Shixiang and Ma, Shiqian and Man-Cho So, Anthony and Zhang, Tong},
title = {Proximal Gradient Method for Nonsmooth Optimization over the Stiefel Manifold},
journal = {SIAM Journal on Optimization},
volume = {30},
number = {1},
pages = {210-239},
year = {2020},
doi = {10.1137/18M122457X},
URL = {https://doi.org/10.1137/18M122457X},
eprint = {https://doi.org/10.1137/18M122457X}
}

@book{Boumal_2023, 
place={Cambridge}, 
title={An Introduction to Optimization on Smooth Manifolds}, 
publisher={Cambridge University Press}, 
author={Boumal, Nicolas}, 
year={2023}}

@article{Chen24,
author = {Chen, Shixiang and Ma, Shiqian and Man-Cho So, Anthony and Zhang, Tong},
title = {Nonsmooth Optimization over the Stiefel Manifold and Beyond: Proximal Gradient Method and Recent Variants},
journal = {SIAM Review},
volume = {66},
number = {2},
pages = {319-352},
year = {2024},
doi = {10.1137/24M1628578},
URL = {https://doi.org/10.1137/24M1628578}}

@article{li2021weakly,
  author = {Li, Xiaoli and Chen, Shixiang and Deng, Zengde and Qu, Qing and Zhu, Zhihui and So, Anthony Man-Cho},
  title = {Weakly Convex Optimization over Stiefel Manifold Using Riemannian Subgradient-Type Methods},
  journal = {SIAM Journal on Optimization},
  volume = {31},
  number = {3},
  pages = {1605--1634},
  year = {2021},
  publisher = {Society for Industrial and Applied Mathematics},
  doi = {10.1137/20M1358517}
}

@article{CambierAbsilSICOMP16,
author = {Cambier, L\'{e}opold and Absil, P.-A.},
title = {Robust Low-Rank Matrix Completion by Riemannian Optimization},
journal = {SIAM Journal on Scientific Computing},
volume = {38},
number = {5},
pages = {S440-S460},
year = {2016},
doi = {10.1137/15M1025153},
URL = {https://doi.org/10.1137/15M1025153},
eprint = {https://doi.org/10.1137/15M1025153}
}

@article{yang2014optimality,
  title={Optimality conditions for the nonlinear programming problems on Riemannian manifolds},
  author={Yang, Wei Hong and Zhang, Lei-Hong and Song, Ruyi},
  journal={Pacific Journal of Optimization},
  volume={10},
  number={2},
  pages={415--434},
  year={2014}
}

@article{XiaPenaltySIOPT21,
author = {Xiao, Nachuan and Liu, Xin and Yuan, Ya-xiang},
title = {Exact Penalty Function for $\ell_{2,1}$ Norm Minimization over the Stiefel Manifold},
journal = {SIAM Journal on Optimization},
volume = {31},
number = {4},
pages = {3097-3126},
year = {2021},
doi = {10.1137/20M1354313},
URL = {https://doi.org/10.1137/20M1354313},
eprint = {https://doi.org/10.1137/20M1354313}
}

@article{Fan_Pall_1957,
title={Imbedding Conditions for Hermitian and Normal Matrices}, 
volume={9}, DOI={10.4153/CJM-1957-036-1}, 
journal={Canadian Journal of Mathematics}, 
author={Fan, Ky and Pall, Gordon}, 
year={1957}, 
pages={298–304}}

@article{AttCzar18,
	author = {Attouch, Hedy and Cabot, Alexandre and Czarnecki, Marc-Olivier},
	journal = {Transactions of the American Mathematical Society},
	number = {2},
	pages = {755--790},
	title = {Asymptotic behavior of nonautonomous monotone and subgradient evolution equations},
	volume = {370},
	year = {2018}}

@article{Stewart_1977,
  author  = {Stewart, G. W.},
  title   = {On the Perturbation of Pseudo-Inverses, Projections and Linear Least Squares Problems},
  journal = {SIAM Review},
  volume  = {19},
  number  = {4},
  pages   = {634--662},
  year    = {1977},
  doi     = {10.1137/1019094}}
\end{document}